\documentclass[onefignum,onetabnum]{siamart190516}

\usepackage{caption}
\usepackage{subcaption}
\usepackage{lipsum}
\usepackage{amsfonts}
\usepackage{graphicx}
\usepackage{epstopdf}
\usepackage{algorithmic}
\ifpdf
  \DeclareGraphicsExtensions{.eps,.pdf,.png,.jpg}
\else
  \DeclareGraphicsExtensions{.eps}
\fi

\newsiamremark{remark}{Remark}
\newsiamremark{hypothesis}{Hypothesis}
\crefname{hypothesis}{Hypothesis}{Hypotheses}
\newsiamthm{claim}{Claim}

\headers{An interior point  method for DF optimization}{A. Brilli, G. Liuzzi, and S. Lucidi}

\title{An interior point  method for nonlinear constrained derivative-free optimization
}

\author{Andrea Brilli\thanks{``Sapienza'' University of Rome, Department of Computer Control and Management Engineering ``A. Ruberti'', Rome, Italy
  (\email{brilli.1708157@studenti.uniroma1.it}, \email{liuzzi@diag.uniroma1.it}, \email{lucidi@diag.uniroma1.it}).}
\and Giampaolo Liuzzi\footnotemark[2]
\and Stefano Lucidi\footnotemark[2]}

\usepackage{amsopn}

\newcommand{\nn}{\nonumber}

\newcommand{\bd}{\begin{document}}
\newcommand{\ed}{\end{document}}

\def\la{\lambda}
\def\eps{\epsilon}
\def\xstar{x^{\star}}
\def\xbar{{\bar x}}
\def\xhat{{\hat x}}
\def\ubar{{\bar u}}
\def\xtilde{{\tilde x}}
\def\utilde{{\tilde u}}
\def\xstar{{x^*}}
\def\ustar{{u^*}}

\renewcommand{\Re}{\mathbb{R}}

\begin{document}

\maketitle

\begin{abstract}
In this paper we consider constrained optimization problems where both the objective and constraint functions are of the black-box type. Furthermore, we assume that the nonlinear inequality constraints are non-relaxable, i.e. their values and that of the objective function cannot be computed outside of the feasible region. This situation happens frequently in practice especially in the black-box setting where function values are typically computed by means of complex simulation programs which may fail to execute if the considered point is outside of the feasible region. For such problems, we propose a new derivative-free optimization method which is based on the use of a merit function that handles inequality constraints by means of a log-barrier approach and equality constraints by means of a quadratic penalty approach. We prove convergence of the proposed method to KKT stationary points of the problem under quite mild assumptions. Furthermore, we also carry out a preliminary numerical experience on standard test problems and comparison with a state-of-the-art solver which shows efficiency of the proposed method.
\end{abstract}



\section{Introduction}
\label{sec:intro}

In this paper we consider the nonlinear constrained minimization
problem
\begin{equation}
\label{prob1}
\begin{array}{l}
\min f(x), \\
\quad g(x)\le 0, \\
\quad h(x)= 0, \\
\quad l\le x\le u, \end{array}
\end{equation}
where $f:\Re^n\to \Re$, $g:\Re^n\to \Re^m$, $h:\Re^n\to\Re^q$, and $l,u\in\Re^n$, with $l<u$,
are vectors of lower and upper bounds on the variables
$x\in\Re^n$. Furthermore, we
assume that $f$, $g$ and $h$ are continuously differentiable functions
even though their derivatives can be neither calculated nor
explicitly approximated. We denote {by ${\cal S}$ the set defined by the nonlinear inequality constraints and} by $X$ the set defined by simple bounds on
the variables, that is,
$$
{\cal S}=\{x\in \Re^n: g(x)\le 0\},
$$
$$
X=\{x\in \Re^n: l\le x\le u\},
$$
and by $\cal F$ the feasible set of problem (\ref{prob1}), namely,
$$
{\cal F} = \{x\in \Re^n: h(x) = 0\}\cap {\cal S}\cap X.
$$
Furthermore, we assume that a point $x_0\in\stackrel{\circ}{\cal S}$ exists.
We note that, by definition, $X$ is a compact set so that $\cal F$ is compact as well.

To solve problem (\ref{prob1}), we resort to  the following merit function in which inequality constraints are handled by log-barrier penalty terms whereas equality constraints are addressed by standard exterior penalty terms (see e.g. \cite{fiaccomccormick}).
$$
P(x;\epsilon) = f(x) - \epsilon\displaystyle{\sum_{j=1}^m} \log\left[-g_j(x)\right]+ \frac{1}{\epsilon} {\sum_{j=1}^q} |h_j(x)|^\nu,
$$
where $\nu>1$. Note that only the nonlinear constraints have been taken into account. Bound constraints on the variables will be addressed explicitly by the optimization algorithm.

We assume that, $P(x;\epsilon) = +\infty$, for all $x\in\Re^n$ such that  $g(x) \not < 0$.
Then, we consider the problem
\cite{bertsekas82,bertsekas99}
\begin{equation}
\label{prob2}
\begin{array}{l}
 \min\ P(x;\eps)\\
 s.t.\ x\in \stackrel{\circ}{\cal S}\cap X
\end{array}
\end{equation}
For every fixed value of the penalty parameter $\eps$, $P(x;\eps)$ is continuously differentiable in $\stackrel{\circ}{\cal S}$ under the stated
assumptions.

In many engineering problems, the values
of the functions defining the objective and constraints of the
problem are computed by means of complex simulation programs. For
this reason, their analytic expressions are not available. Hence, derivatives are not available or, at the very least, they are untrustworthy. 

Many real world applications fit into the derivative-free or black-box optimization paradigm. Such problems usually present nonlinear constraints along with bound constraints on the variables. Black-box optimization problems are widely studied in the literature (see, e.g., \cite{audet2017derivative,conn2009introduction,torczon:03}) and many algorithms have been proposed for the solution of constrained black-box optimization problems. 
In particular, in \cite{lewis2002globally} the use of an augmented Lagrangian function in connection  with a pattern search algorithm has been proposed. In \cite{liuzzi2010sequential} a sequential penalty derivative-free linesearch approach has been studied, whereas in \cite{liuzzi2009derivative} the use of a nonsmooth exact penalty function has been proposed. A mesh adaptive direct search method, namely NOMAD, has been firstly introduced and analyzed in \cite{audet2006mesh} to solve constrained black-box problems by using an extreme penalty function to manage general and hidden constraints. 

According to \cite{larson:2019,digabel2015taxonomy} inequality constraints can be either relaxable or unrelaxable. Unrelaxable constraints are those constraints that must always be {satisfied} by the points produced by the optimization algorithm. Hence, when unrelaxable balck-box constraints are present, the optimization algorithm should take into proper account this feature. Typically, such constraints can be managed by a so-called {\em extreme} or {\em death} penalty approach (see e.g. \cite{NOMAD}). In particular, an objective function value of $+\infty$ is assigned to points that are unfeasible with respect to one or more unrelaxable constraints. However, it should also be mentioned that such penalization strategy, by making the objective function discontinuous on the boundary of the feasible region, introduces many difficulties and ill-conditioning in the problem. As a result, solving the problem could become impractical or, at the very least, the computed solution could be far away from the real solution point. 

A possible way of handling the above mentioned difficulty, consists in the use of some sort of interior penalization that modifies the landscape of the objective function in the interior of the feasible region by adding to the objective function terms that gradually tend to $+\infty$ as the points approach the boundary of the feasible region (see e.g. \cite{forsgreen:interior,curtis2012penalty,nocedal2006numerical}).

\par\smallskip

The paper is organized as follows. In Section \ref{sec:notation}, we introduce some notation and preliminary results that will be used in the paper. Section \ref{sec:3} is devoted to the definition of a minimization algorithm for the proposed merit function when the barrier parameter is held fixed. Also, quite standard convergence analysis is reported. In section \ref{sec:4}, the main algorithm is described along with its theoretical convergence analysis. 
Section \ref{sec:experiments} is devoted to the numerical experimentation and comparison of the proposed method with a stat-of-the-art solver, namely NOMAD \cite{NOMAD}. In Section \ref{sec:conclusion} we draw some conclusions. Finally, in Appendix \ref{app:tech}, more technical results (which are used to prove convergence of the main algorithm) are proved.

\section{Notation and preliminary results}
\label{sec:notation}

In this section we introduce some notation and assumptions
that will be used throughout the paper.

Given a vector $v\in\Re^n$, a subscript will be used to denote
either one of its components ($v_i$) or the fact that it is an
element of an infinite sequence of vectors ($v_k$). To avoid
possible misunderstanding or ambiguities, the $i$th component of a
vector will be denoted by $(v)_i$. We denote by $v^j$ the generic
$j$th element of a finite set of vectors. Given two vectors
$a,b\in\Re^n$, we denote by $y=\max\{a,b\}$ the vector such that
$y_i=\max\{a_i,b_i\}$, $i=1,\dots,n$. Furthermore, given a vector
$v$, we denote $v^+=\max\{0,v\}$.

\begin{definition}[cone of feasible directions]
\label{cone}
Given a point $x\in X$, let
$$
D(x)=\{d\in \Re^n:d_i\ge 0~\mbox{if}\ x_i = l_i,\ d_i\le 0~\mbox{if}\ x_i = u_i,\ i=1,\dots,n\}
$$
be the cone of feasible directions at $x$ with respect to the simple bound constraints.
\end{definition}

Let $L(x,\la,\mu)$ be
the Lagrangian function associated with the nonlinear constraints of problem (\ref{prob1}),
$$
 L(x,\lambda,\mu)= f(x)+\lambda^Tg(x)+\mu^Th(x)
$$

We recall the Mangasarian--Fromovitz constraint qualification
(MFCQ).

\begin{definition}
\label{ass}
A point $x\in X$ is said to satisfy the MFCQ if two conditions are satisfied:
\begin{itemize}
    \item[\rm (a)] There does not exist a nonzero vector $\alpha=(\alpha_1,...,\alpha_q)$ such that:
    \begin{equation}
    \label{box_mfcq} \left(\sum_{i=1}^{q} \alpha_i\nabla h_i(x)\right)^T d \geq 0,\quad\quad\forall d\in D(x),
    \end{equation}
    
    \item[\rm (b)] there exists a feasible direction $d\in D(x)$:
    \begin{equation}
    \nabla g_l(x)^T d < 0\ \ \ \forall l\in I(x),\quad
    \nabla h_j(x)^Td=0\ \ \ \forall j=1,...,q
    \end{equation}
where $I(x)=\{i:g_i(x)\geq0\}$.
\end{itemize}
\end{definition}

The following proposition is a well-known result (see, for
instance, \cite{bertsekas99}) which states necessary optimality
conditions for problem (\ref{prob1}).

\begin{proposition}
\label{bert} Let $x^\star\in\cal F$ be a local minimum of problem
{\rm (\ref{prob1})} that satisfies the MFCQ. Then, there exists a vectors $\lambda^\star\in
\Re^m$, $\mu^\star\in\Re^q$ such that
\begin{equation}
\label{derdir} \nabla_x L(x^\star,\lambda^\star,\mu^\star)^T (x-x^\star)\geq
0\quad\quad\forall x\in X,
\end{equation}
\begin{equation}
\label{cmp} (\lambda^\star)^Tg(x^\star)=0,\quad\quad
\lambda^\star\ge 0.\qquad \hbox{{\vbox{\hrule height0.6pt\hbox{%
   \vrule height1.3ex width0.6pt\hskip0.8ex
   \vrule width0.6pt}\hrule height0.6pt
  }}}
\end{equation}
\end{proposition}
\unskip

\begin{definition}[{stationary point}]
\label{stationary} A point $x^\star\in\cal F$ is said to be a
stationary point for problem $(\ref{prob1})$ if a vector
$\lambda^\star\in R^m$ and $\mu^\star\in\Re^q$ exists such that $(\ref{derdir})$ and
$(\ref{cmp})$ are satisfied.
\end{definition}

{ Now  we recall two  results from \cite{LinLucidiPalagiSciandrone} and \cite{liuzzi2010sequential}
concerning the set $D(x)$.}

\begin{proposition}
\label{iota} Let $\{x_k\}$ be a sequence of points such that
$x_k\in X$ for all $k$. Assume further that $x_k \to \bar x$ for
$k\to\infty$. Then, given any direction $\bar d\in D(\bar x)$,
there exists a scalar $\bar\beta > 0$ such that, for sufficiently
large $k$, we have
\[
x_k+\beta\bar d\in X\quad \forall\beta\in [0,\bar\beta].
\]
\end{proposition}

Hence, given a sequence $\{x_k\}\subset X$ such that $x_k\to\bar x$ for $k\to\infty$, it results $D(\bar x)\subseteq D(x_k)$ for $k$ sufficiently large.




Now we define the set of unit vectors
$$
D = \{\pm e^1,\dots,\pm e^n\},
$$
where $e^i$, $i=1,\dots,n$, is the $i$th unit coordinate vector. 

{The following proposition  shows} that set
$D$ contains the generators of the cone of feasible directions $D(x)$ at any point $x\in X$.

\begin{proposition}\label{set D}
\label{cn1}
Let $x\in X$. We have
\begin{equation}
\label{spdir}
cone\{D\cap D(x)\}=D(x).
\end{equation}
\end{proposition}

\section{Minimization of $P(x;\eps)$ when $\eps$ is fixed} 
\label{sec:3}


In this section we define and analyse an almost ``classical'' derivative-free algorithm with linesearches for the solution of problem (\ref{prob2}) when the penalty parameter $\eps$ is kept fixed. Such an algorithm is reported in the box below.
\par\medskip

\noindent\framebox[\textwidth]{\parbox{0.95\textwidth}{\small
\par
\centerline{ {\bf Derivative-free linesearch method (DFL)}.}
\par\medskip

 {\bf Data.} $x_0\in X$ such that $g(x_0)<0$, $\epsilon>0$, $\gamma > 0$,
 $\theta\in (0,1)$, $p>1$, $\tilde\alpha_0^i > 0$,

 \par\quad\quad\quad and set $d_0^i=e^i$ for $i=1,\ldots,n$.

\par\medskip

 {\bf For} $k = 0,1,2,\dots$ {\bf do} {\em (Main iteration loop)}
 \begin{itemize}\item[] {\bf Step 1.} Set $y_k^1=x_k$.
 \begin{itemize}
 \item[ ]{\bf For} $i = 1,\dots,n$ {\bf do} {\em (Exploration of the search directions)}
 \begin{itemize}
 \item[ ]{\bf Step 1.2.} Compute 
 $\hat\alpha_k^i\leq \tilde\alpha_{k}^i$ s.t. $y_k^i+\hat\alpha_k^i d_k^i\in\stackrel{\circ}{\cal S}\cap X$ 
\par\quad {\bf If} $\hat\alpha_k^i>0$,  and $P(y_k^i+\hat\alpha_k^i d_k^i;\epsilon)\le
P(y_k^i;\epsilon)-\gamma{(\hat\alpha_k^i)}^2$,
 \par\qquad compute $\alpha_k^i$ by the {\it Expansion
 Step($\hat\alpha_k^i,y_k^i,d_k^i,\gamma;\alpha_k^i)$};
 \par\qquad set
 $\tilde\alpha_{k+1}^{i}=\alpha_k^i$, $d^i_{k+1}=d^i_{k}$ and go to {\bf Step 1.5.}

\item[ ]{\bf Step 1.3.} Compute 
 $\hat\alpha_k^i\leq \tilde\alpha_{k}^i$ s.t. $y_k^i-\hat\alpha_k^i d_k^i\in\stackrel{\circ}{\cal S}\cap X$ 
\par\quad {\bf If} $\hat\alpha_k^i>0$,  and $P(y_k^i-\hat\alpha_k^i d_k^i;\epsilon)\le
P(y_k^i;\epsilon)-\gamma{(\hat\alpha_k^i)}^2$,
 \par\qquad compute $\alpha_k^i$ by the {\it Expansion
 Step($\hat\alpha_k^i,y_k^i,-d_k^i,\gamma;\alpha_k^i)$};

 \par\qquad set
 $\tilde\alpha_{k+1}^{i}=\alpha_k^i$, $d^i_{k+1}=-d^i_{k}$, and go to {\bf Step 1.5.}

 \item[ ] {\bf Step 1.4.} Set $\alpha_k^i=0$, $d_{k+1}^i = d_k^i$, $\tilde\alpha_{k+1}^{i}=\theta\tilde\alpha_k^i$

 \item[ ] {\bf Step 1.5.} Set $y_k^{i+1}=y_k^i+\alpha_k^id_k^i$. 
 \end{itemize}
 \item[ ] {\bf Endfor} 
 \end{itemize}

 \item[] {\bf Step 3.} Find $x_{k+1}\in \stackrel{\circ}{\cal S}\cap X$
 such that $P(x_{k+1};\epsilon)\le P(y_k^{n+1};\epsilon)$.
 \end{itemize}
 {\bf Endfor}
\par
\noindent }}
\par\bigskip

The minimization process of such a derivative-free method is based on  suitable sampling techniques
along a set of directions that are able to convey, in the limit,
sufficient knowledge of the problem functions to recover first
order information. In particular, for  box constrained  optimization problems,  suitable choice for set of directions are the unit coordinate vectors $e^i$, $i=1,\dots,n$, (see {Proposition \ref{set D} and } \cite{torczon:03,LucidiSciandrone:02}). Indeed, the search directions are  initialized to the unit coordinate vectors, i.e. $d_0^i = e^i$, for $i=1,\dots,n$. Then, at iteration $k$, the algorithm defines the directions to be used at iteration $k+1$. More in particular, for $i\in\{1,\dots,n\}$, the following rule is adopted to define $d_{k+1}^i$: 
\[
    d_{k+1}^i = \left\{\begin{array}{ll}
       \phantom{-}d_k^i & \mbox{when ``sufficient'' decrease is achieved along $d_k^i$,}\\
       -d_k^i & \mbox{when ``sufficient'' decrease is achieved along $-d_k^i$,}\\
       \phantom{-}d_k^i & \mbox{otherwise.}
    \end{array}\right.
\]

It is worth noticing that the sole source of complexity in the proposed algorithm (with respect to standard linesearch-based derivative-free algorithms for bound constrained problems) resides in the {\em Expansion Step}, where appropriate actions must be performed to take into account that the objective function can be computed only on $\stackrel{\circ}{\cal S}\cap X$. The Expansion Step procedure is reported below.

\par\medskip

\noindent\framebox[\textwidth]{\parbox{0.95\textwidth}{\small
\par
\centerline{{\bf Expansion Step
(\boldmath$\hat\alpha,y,p,\gamma;\alpha$).}}
\begin{itemize}
 \item[ ] {\bf Data.} $\delta\in (0,1)$ and $b$ the largest step such that $y+bp\in X$.
 \item[ ] {\bf Step 1.} \rm{Set $\alpha\leftarrow\hat\alpha$.
 \item[ ] {\bf Step 2.} set $\check\alpha\leftarrow\min\{b,{\alpha}/{\delta}\}$
 \item[ ] {\bf Step 3.} {\bf If} $y+\check\alpha p\not\in\stackrel{\circ}{\cal S}$ {\bf return}
 \item[ ] \phantom{\bf Step 3.}  {\bf Elseif} $\check \alpha < b$ and
 $P(y+\check\alpha p;\epsilon)\leq P(y;\epsilon)-\gamma\check\alpha^2$ {\bf then}  
 \item[ ] \phantom{\bf Step 3.} \qquad  set $\alpha \leftarrow \check\alpha$
 \item[ ] \phantom{\bf Step 3.}  {\bf Elseif} $\check \alpha = b$ and  
 $P(y+\check\alpha p;\epsilon) \le P(y;\epsilon)-\gamma\check\alpha^2$ {\bf then}
 \item[ ] \phantom{\bf Step 3.} \qquad  set $\alpha \leftarrow \check\alpha$ and {\bf return}
 \item[ ] \phantom{\bf Step 3.}  {\bf Else} (i.e. $P(y+\check\alpha p;\epsilon)> P(y;\epsilon)-\gamma\check\alpha^2$)
 \item[ ] \phantom{\bf Step 2.} \qquad {\bf return}
 \item[ ] {\bf Step 4.} Go to Step 2.
 }
\end{itemize}
\par
\noindent }}
\par\bigskip\noindent

As we can see, the {\em Expansion Step} procedure is invoked when it is possible to find a strictly positive initial stepsize $\hat\alpha$ which gives sufficient decrease, i.e. $\hat\alpha > 0$ and
\[
 P(y+\hat\alpha p;\eps) \leq P(y;\eps) -\gamma\hat\alpha^2
\]
Then, the {\em Expansion Step} procedure computes two stepsizes, namely $\alpha$ and $\check\alpha$. The stepsize $\alpha$ is such that $\alpha \geq \hat\alpha$ and sufficient decrease is obtained w.r.t. the initial point, i.e. 
\[
 P(y+\alpha p;\eps) \leq P(y;\eps)-\gamma\alpha^2.
\]
Furthermore, the Expansion Step also defines a stepsize $\check\alpha$ such that:
\begin{itemize}
    \item[i)] either $y+\check\alpha p\not\in\stackrel{\circ}{\cal S}$;
    \item[ii)] or $\check\alpha=b$, $y+\check\alpha p\in X\cap \stackrel{\circ}{\cal S}$ and
    \[
    P(y+\check\alpha p;\epsilon) \le P(y;\epsilon)-\gamma\check\alpha^2;
    \]
    \item[iii)] or $\check\alpha = \min\{b,\alpha/\delta\}$,  $y+\check\alpha p\in X\cap \stackrel{\circ}{\cal S}$ and
    \[
    P(y+\check\alpha p;\epsilon) > P(y;\epsilon)-\gamma\check\alpha^2.
    \]    
\end{itemize}
Note that, point (i) implies that a stepsize $\beta < \check\alpha$ exists (see the proof of Proposition \ref{part_a_theorem}) such that  $y+\beta p\in X\cap \stackrel{\circ}{\cal S}$ and
    \[
    P(y+\beta p;\epsilon) > P(y;\epsilon)-\gamma\beta^2.
    \] 
Point (ii), eventually, cannot happen (see the proof of Proposition \ref{part_a_theorem}). Point (iii) is the interesting case and the one that will be used in the proof of Proposition \ref{part_a_theorem}.

\subsection{Convergence analysis}

In this subsection we carry out the convergence analysis of Algorithm DFL (e.g. the minimization algorithm for $P(x;\eps)$ when parameter $\eps$ is fixed). In the following proposition, we first prove that the Expansion Step procedure is well defined.
\begin{proposition}
The {\em Expansion Step} is well defined, i.e. it always returns a step $\alpha$.
\end{proposition}
\begin{proof}
We proceed by contradiction and assume that the procedure infinitely cycles. If that is the case, it produces an infinite sequence of values \{$\check\alpha_j$\}. By the instructions:
\[
\check\alpha_j = \frac{\alpha}{\delta^j},
\]
which contradicts $\check\alpha < b$.
\end{proof}
Now, we prove that the sequences of tentative stepsizes $\{\tilde\alpha_k^i\}$ and actual stepsizes $\{\alpha_k^i\}$, for $i=1,\dots,n$, are all convergent to zero.

\begin{proposition}
\label{eps_to_0_1}
Let $\{\alpha_k^i\}$ and $\{\tilde\alpha_k^i\}$, $i=1,\dots,n$, be the sequences produced by the Algorithm, then
\begin{equation}
\label{utile1n_new_1} \lim_{k\to\infty}\alpha_k^i=0\quad\quad {\it
for~}i=1,\ldots,n,
\end{equation}
\begin{equation}
\label{utile0n_new_1} \lim_{k\to\infty}\tilde\alpha_k^i=0\quad\quad {\it
for~}i=1,\ldots,n,
\end{equation}
\end{proposition}

\begin{proof}
For every $i=1,\ldots,n$ we prove (\ref{utile1n_new_1}) by splitting the
iteration sequence $\{k\}$ into two parts, $K^{\prime}$ and $K^{\prime\prime}$.
We identify with $K^{\prime}$ those iterations where
\begin{equation}
\label{zero_new_1}
\alpha_k^i=0,
\end{equation}
and with $K^{\prime\prime}$ those iterations where $\alpha_k^i\neq 0$ is
produced by the Expansion Step. Then the instructions of the
algorithm imply
\begin{eqnarray}
\label{dimin_new_1} \\
P(x_{k+1};\eps)\le P(y_k^i+\alpha_k^i d_k^i;\eps)\leq P(y_k^i;\eps)-\gamma{(\alpha_k^i)}^2{\|d_k^i\|}^2\leq
P(x_k;\eps)-\gamma{(\alpha_k^i)}^2{\|d_k^i\|}^2. \nn
\end{eqnarray}
Taking into account the compactness assumption on $X$, it follows
from (\ref{dimin_new_1}) that $\{P(x_k;\eps)\}$ tends to a limit
$\bar P$. If $K^\prime$ is infinite, then from (\ref{zero_new_1}) we trivially have that
\[
\lim_{k\to\infty, k\in K^\prime} \alpha_k^i =0
\]
If, on the other hand, $K^{\prime\prime}$ is an infinite subset, recalling that
$\|d_k^i\|=1$, we obtain
\begin{equation}
\label{onalfapla_new_1}
\lim_{k\to\infty,k\in K^{\prime\prime}}\alpha_k^i=0.
\end{equation}
Therefore, (\ref{zero_new_1}) and (\ref{onalfapla_new_1}) imply (\ref{utile1n_new_1}).\par
In order to prove (\ref{utile0n_new_1}), for each
$i\in\{1,\ldots,n\}$ we split the iteration sequence $\{k\}$ into
two parts, $K_1$ and $K_2$. We identify with $K_1$ those
iterations where the Expansion Step has been performed using the
direction $d_k^i$, for which we have
\begin{equation}
\label{a1_new_1}
\tilde\alpha_{k+1}^i=\alpha_k^i.
\end{equation}
We denote by $K_2$ those iterations
where we have failed in decreasing the objective function
 along the directions $d_k^i$ and $-d_k^i$. By the instructions
of the algorithm it follows that for all $k\in K_2$
\begin{equation}
\label{diminuzione_new_1}
\tilde\alpha_{k+1}^i\le\theta\tilde\alpha_k^i,
\end{equation}
where $\theta\in (0,1)$.\par
If $K_1$ is an infinite
subset, from (\ref{a1_new_1}) and (\ref{utile1n_new_1}) we get that
\begin{equation}
\label{onalfatildepla_new_1}
\lim_{k\to\infty,k\in K_1}\tilde\alpha_{k+1}^i=0.
\end{equation}
Now, let us assume that $K_2$ is an infinite subset. For each
$k\in K_2$, let $m_k$ (we omit the dependence on $i$) be the
biggest index such that $m_k<k$ and $m_k\in K_1$. Then we have
\begin{equation}
\label{minorpla_new_1}
\tilde\alpha_{k+1}^i\le
\theta^{(k+1-m_k)}\tilde\alpha_{m_k}^i
\end{equation}
(we can assume $m_k=0$ if the index $m_k$ does not exist, that is,
$K_1$ is empty).\par
As $k\to\infty$ and $k\in K_2$, either $m_k\to\infty$ (namely,
$K_1$ is an infinite subset) or $(k+1-m_k)\to\infty$ (namely,
$K_1$ is finite). Hence, if $K_2$ is an infinite subset,
(\ref{minorpla_new_1}) together with (\ref{onalfatildepla_new_1}), or the fact
that $\theta\in (0,1)$, yields
\begin{equation}
\label{ovviopla_1}
\lim_{k\to\infty,k\in K_2}\tilde\alpha_{k+1}^i=0,
\end{equation}
so that (\ref{utile0n_new_1}) is proved.
\end{proof}

Then, we report a technical proposition that states general convergence conditions that will be used in the proof of the main convergence theorem.

\begin{proposition}\label{part_a_theorem_1}
Let $\{x_k\}$ and $\{y_k^i\}$, $i=1,\dots,n+1$, be the sequences produced by Algorithm DFL {and let $\{x_k\}_{H}$ be a subsequence converging to the point $\bar x$. Then,  for $k\in {H}$} sufficiently large, for all $d^i\in D\cap D(\bar x)$,
there exist 
scalars $\xi_k^i>0$ such that
\begin{eqnarray}
 \label{d0bis_old_1}
 && y_k^i+ \xi_k^id^i\in X\cap \stackrel{\circ}{\cal S},\\
 %
 \label{d1bis1_old_1}
 && P({y_k^i}+ \xi_k^i d^i;\eps)\ge P({y_k^i};\eps)-o(\xi_k^i),\\
 %
 \label{xi_ki_to_zero_1}
    &&  \lim_{k\to\infty{, k\in{H}}} \xi_k^i = 0\\
\label{y_ki_to_xbar_1}
    && \lim_{k\to\infty{, k\in{H}}} \|y_k^i - x_k\| = 0
\end{eqnarray}
 
\end{proposition}
\begin{proof}
The proof of this proposition is reported in the next section (see Proposition \ref{part_a_theorem}) for the more general setting when a sequence $\{\eps_k\}$ of penalty-barrier parameters is considered.
\end{proof}

Finally, we report the main convergence result concerning the linesearch algorithm for the solution of problem (\ref{prob2}).

\begin{theorem}
Let $\{x_k\}$ be the sequence of points produced by the algorithm. Then, every limit point $\bar x$ of $\{x_k\}$ is stationary for problem (\ref{prob2}), namely
\[
  \nabla P(\bar x;\eps)^\top (x-\bar x) \geq 0,\quad\forall\ x\in \stackrel{\circ}{\cal S}\cap X.
\]
\end{theorem}
\begin{proof}
Since $\{x_k\}\subseteq X$ and $X$ is compact, the sequence $\{x_k\}$ admits limit points. Let us consider one such limit point $\bar x$, i.e. an index set $\bar K$ exists such that 
\[
\lim_{k\to\infty,k\in \bar K} x_k = \bar x.
\]
Let us denote $\bar D=D\cap D(\bar x)$. {Recalling Proposition \ref{part_a_theorem_1}
we have that (\ref{d0bis_old_1}), (\ref{d1bis1_old_1}),   (\ref{xi_ki_to_zero_1}) and (\ref{y_ki_to_xbar_1}) hold.}

\par\medskip\noindent

Recalling (\ref{dimin_new_1}) and taking into account the compactness assumption on X, it follows that $\{P(x_k;\eps)\}$ tends to a limit $\bar P$. Since
\[
\lim_{k\to \infty, k\in\bar K}g_\ell(x_k)=0\ \mbox{for some}\ \ell\in\{1,\dots,m\} \quad \text{implies } \lim_{k\to \infty, k\in\bar K}P(x_k;\eps) = +\infty,
\]
we have that $g_\ell(\bar x)<0,$ for all $\ell=1,\dots,m$. Thus, an $r>0$ exists  such that $g_\ell(x)<0$ for all $x\in B(\bar x;r)$ and $\ell=1,\dots,m$. Recalling (\ref{xi_ki_to_zero_1}) and (\ref{y_ki_to_xbar_1}) we can state that an index $\bar k$ exists such that for all $k\ge\bar k,k\in \bar K$:
\[
y_k^i+\xi_k^id^i \in B(\bar x; r),
\]
which implies:
\[
g_\ell(y_k+t\xi_k^id^i)<0 \ \forall t\in[0,1],\forall k\ge\bar k, k\in\bar K, \forall \ell=1,\dots,m.
\]

\noindent
We have that for all $k\in\bar K$ sufficiently large and for all $d^i\in D\cap{D}(\bar{x})$:
\begin{itemize}
    \item[i)] $y_k^i\in X\cap \stackrel{\circ}{\cal S}$
    \item[ii)]$y_k^i+\xi_k^id^i\in X\cap \stackrel{\circ}{\cal S}$
\end{itemize}
Now, let us assume there exists $\hat{t}_k^i \in (0,1)$ such that $y_k^i+\hat{t}_k^i \xi_k^i d^i \notin  X\cap \stackrel{\circ}{\cal S}$. Then, by the compactness of $X$
\[
y_k^i+\hat{t}_k^i\xi_k^id^i \notin \stackrel{\circ}{\cal S}, \quad \text{i.e. } g_l(y_k^i+\hat{t}_k^i\xi_k^id^i)>0
\]
Using the continuity assumption on the constraints, there exists at least one constant $\check{t}_k^i \in (0,\hat{t}_k^i)$ such that:
\[
g_l(y_k^i+\check{t}_k^i \xi_k^id^i)=0,
\]
if multiple constants exist that satisfy the condition above, we will consider $\check{t}_k^i$ to be the smallest one. We have now that $y_k^i + t\xi_k^id^i \in X\cap \stackrel{\circ}{\cal S}$ for all $t \in [0,\check{t}_k^i)$, so that by the definition of $P$:
\[
P(y;\eps_k) \text{ is continuous } \forall y \in [y_k^i,y_k^i+\check{t}_k^i\xi_k^id^i) \text{ and } P(y_k^i+\check{t}_k^i\xi_k^id^i;\eps_k)=+\infty,
\]
thus, a constant $t_k^{*,i} \in (0,\check{t}_k^i)$ must exist such that:
\begin{eqnarray}
    && \label{new_fallimento1}P({y_k^i}+ t_k^{*,i} \xi_k^i d^i;\eps_k)\ge P({y_k^i};\eps_k)-o(t_k^{*,i}\xi_k^i)\\
    && \label{new_feasibility1}y_k^i+t \xi_k^i d^i \in X\cap \stackrel{\circ}{\cal S} \forall t\in[0,t_k^{*,i}].
\end{eqnarray}
Let us now set
\[
\bar{\xi_k^i} = t_k^{*,i} \xi_k^i \le \xi_k^i
\]
and apply the mean-value theorem to (\ref{new_fallimento1}). Thus, we can write

$$
-o(\bar\xi_k^{i})\le P(y_k^i+\bar\xi_k^{i} d^{i};\eps_k) -
P(y_k^i;\eps_k)= \bar\xi_k^{i}\nabla
P(u_k^i;\eps_k)^{T}d^{i}\quad\quad \forall d^i\in \bar D,
$$
where $u_k^i=y_k^i+t_k^{i}\bar\xi_k^{i} d^{i}$, with $t_k^{i}\in (0,1)$. 
Thus, we have
\[
-\frac{o(\bar\xi_k^{i})}{\bar\xi_k^{i}} \leq \nabla P(u_k^i;\eps_k)^{T}
d^i\quad\quad \forall d^i\in \bar D.
\]
Recalling Proposition \ref{set D}, (\ref{xi_ki_to_zero_1}) and (\ref{y_ki_to_xbar_1}), taking the limit for $k\to\infty, k\in \bar K$, we get
\[
\nabla P(\bar x;\eps)^\top d \geq 0,\quad \forall\ d\in D(\bar x).
\]
\end{proof}

\section{The main algorithm}\label{sec:4}
In this section, we report the scheme of the proposed derivative-free algorithm used to solve the constrained problem (\ref{prob1}) by means of the log-barrier penalty function $P(x;\eps)$. The scheme of the algorithm is obtained by suitably modifying the algorithm introduced in the previous section which solves problem (\ref{prob2}) instead. In particular, as expected, the only difference consists in the penalty-barier parameter updating rule. Indeed, the following LOG-DFL algorithm exactly is Algorithm DFL except for Step 2 which we highlighted by a box.   



\noindent\framebox[\textwidth]{\parbox{0.95\textwidth}{\small
\par
\centerline{ {\bf Algorithm LOG-DFL}.}
\par\medskip

 {\bf Data.} $x_0\in X$ such that $g(x_0)<0$, $\epsilon_0>0$, $\gamma > 0$,
 $\theta\in (0,1)$, $p>1$, $\tilde\alpha_0^i > 0$,

 \par\quad\quad\quad and set $d_0^i=e^i$ for $i=1,\ldots,n$.

\par\medskip

 {\bf For} $k = 0,1,2,\dots$ {\bf do} {\em (Main iteration loop)}
 \begin{itemize}\item[] {\bf Step 1.} Set $y_k^1=x_k$
 \begin{itemize}
 \item[ ]{\bf For} $i = 1,\dots,n$ {\bf do} {\em (Exploration of the search directions)}
 \begin{itemize}
 \item[ ]{\bf Step 1.2.} Compute 
 $\hat\alpha_k^i\leq \tilde\alpha_{k}^i$ s.t. $y_k^i+\hat\alpha_k^i d_k^i\in\stackrel{\circ}{\cal S}\cap X$ 
\par\quad {\bf If} $\hat\alpha_k^i>0$,  and $P(y_k^i+\hat\alpha_k^i d_k^i;\epsilon_k)\le
P(y_k^i;\epsilon_k)-\gamma{(\hat\alpha_k^i)}^2$,
 \par\qquad compute $\alpha_k^i$ by the
 \par\qquad\quad{\it Expansion
 Step($\hat\alpha_k^i,y_k^i,d_k^i,\gamma$ $;\alpha_k^i$)};
 \par\qquad set
 $\tilde\alpha_{k+1}^{i}=\alpha_k^i$, $d^i_{k+1}=d^i_{k}$ and go to {\bf Step 1.5.}
 
\item[ ]{\bf Step 1.3.} Compute 
 $\hat\alpha_k^i\leq \tilde\alpha_{k}^i$ s.t. $y_k^i-\hat\alpha_k^i d_k^i\in\stackrel{\circ}{\cal S}\cap X$ 
\par\quad {\bf If} $\hat\alpha_k^i>0$,  and $P(y_k^i-\hat\alpha_k^i d_k^i;\epsilon_k)\le
P(y_k^i;\epsilon_k)-\gamma{(\hat\alpha_k^i)}^2$,
 \par\qquad compute $\alpha_k^i$ by the
 \par\qquad\quad {\it Expansion
 Step($\hat\alpha_k^i,y_k^i,-d_k^i,\gamma,$ $;\alpha_k^i$)};
 \par\qquad set
 $\tilde\alpha_{k+1}^{i}=\alpha_k^i$, $d^i_{k+1}=-d^i_{k}$, and go to {\bf Step 1.5.}

 \item[ ] {\bf Step 1.4.} Set $\alpha_k^i=0$, $d_{k+1}^i = d_k^i$, $\tilde\alpha_{k+1}^{i}=\theta\tilde\alpha_k^i$. 
 
 
 \item[ ] {\bf Step 1.5.} Set $y_k^{i+1}=y_k^i+\alpha_k^id_k^i$. 
 \end{itemize}
 \item[ ] {\bf Endfor} 
 \end{itemize}
 
 \fbox{
 \begin{minipage}[t]{0.65\textwidth}
 
 \begin{itemize}
 \item[] \hspace*{-1.5cm}{\bf Step 2.} Set $(g_{min})_k = \displaystyle\min_{i=1,\dots,n+1,\ell=1,\dots,m}\{|g_\ell(y_k^i)|\}$
 \par\medskip
 \hspace*{-1.5cm}\phantom{{\bf Step 2.}} {\bf If} $\max_{i=1,2,\dots,n}\{\tilde\alpha_k^i,\alpha_k^i\}\le\min\{
 \epsilon_k^p$,  
 $(g_{\min})_k^2\}$
 \par\medskip
 \hspace*{-1.5cm}\phantom{{\bf Step 2.}} {\bf Then} choose $\epsilon_{k+1} = \theta\epsilon_k$ {\bf Else} set $\epsilon_{k+1} =
 \epsilon_k$.
 \end{itemize}
\end{minipage}
 }

 \item[] {\bf Step 3.} Find $x_{k+1}\in \stackrel{\circ}{\cal S}\cap X$
 such that $P(x_{k+1};\epsilon_k)\le P(y_k^{n+1};\epsilon_k)$.
 \end{itemize}
 {\bf Endfor}
\par
\noindent }}
\par\bigskip

In particular, about algorithm LOG-DFL, it is worth noting that two quantities are computed during the inner for loop of the algorithm, namely 
\begin{itemize}
    \item[i)] a ``maximum stepsize'' $\max_{i=1,2,\dots,n}\{\tilde\alpha_k^i,\alpha_k^i\}$, i.e. the maximum stepsize used by the algorithm in the entire inner for loop. We recall that this quantity can be roughly considered as a measure of stationarity for the penalty function, see e.g. \cite{LucidiSciandrone:02,torczon:03};
    
    \item[ii)] a ``minimum value'' for the non-relaxable  inequality constraints $(g_{\min})_k$, i.e. the smallest absolute value of the inequality constraints found in the inner loop, namely: $\min_{i=1,\dots,n+1,\ell=1,\dots,m}\{|g_\ell(y_k^i)|\}$.
\end{itemize}
These two quantities play a crucial role in the penalty-barrier parameter updating rule that we shall describe below.

At the end of the inner for loop, the algorithm checks whether the penalty-barrier parameter should be updated. Finally, the new point $x_{k+1}$ is computed by selecting any point which is better than the one produced by the inner for loop.

As concerns the updating rule performed at step 2 of the algorithm, a few comments are in order to help better understand its meaning. 
The algorithm updates the penalty-barrier parameter when the measure of stationarity $\max_{i=1,2,\dots,n}\{\tilde\alpha_k^i,\alpha_k^i\}$ is smaller than the smallest value between $\eps_k^p$ and $(g_{\min})_k^2$. In more details, $\eps_k$ is reduced when both the following conditions are satisfied.  
\begin{itemize}
    \item[i)] $\max_{i=1,2,\dots,n}\{\tilde\alpha_k^i,\alpha_k^i\}$ is  smaller than $\eps_k^p$; 
    \item[ii)] $\max_{i=1,2,\dots,n}\{\tilde\alpha_k^i,\alpha_k^i\}$ is smaller than $(g_{\min})_k^2$.
\end{itemize}
Condition (i) requires that the measure of stationarity is better than the quality of the approximation performed by the merit function w.r.t. 
the constrained problem. \\
On the other hand, condition (ii) requires that the step size used by the algorithm is sufficiently small in order to drive the iterates toward the boundary of the feasible region. \\
It's worth noticing that both conditions imply that the maximum stepsize\\ $\max_{i=1,2,\dots,n}\{\tilde\alpha_k^i,\alpha_k^i\}$ must go to zero faster than the penalty-barrier parameter (which is required to go to zero in order for the iterate to approach a KKT point in the limit) and than the minimum value for the non-relaxable inequality constraints (in the case the limit point lies on the boundary of the feasible region) respectively.


\subsection{Convergence analysis}\label{sec:convanalysis}
This section is devoted to the analysis of the convergence properties of the proposed algorithm.

The next proposition ensures that the updating rule of the algorithm produces a sequence of values of the penalty parameter  which tends  to zero. This result is of paramount importance since the parameter $\eps$ multiplies the log-barrier terms of the merit function.

\begin{proposition}
\label{eps_to_0}
Let $\{\epsilon_k\}$ be the sequence produced by Algorithm LOG-DFL, then
\[
\lim_{k\to\infty} \epsilon_k = 0
\]
\end{proposition}

\begin{proof}
By the instructions of the algorithm, $\{\epsilon_k\}$ is a monotonically non-in\-crea\-sing sequence of positive numbers. Hence, it is convergent to a limit $\bar\eps \geq 0$.
Then, we proceed by contradiction and assume that $\bar\epsilon > 0$. This means that, for $k$ sufficiently large, $\epsilon_k$ is no longer updated. Hence, we can assume that $\eps_k$ stays fixed, i.e. $\eps_k=\bar\eps$, definitely, i.e. the test at step 2 of Algorithm LOG-DFL is no longer satisfied that is
\begin{equation}
\label{utile3n_new}
\max_{i=1,2,\dots,n}\{\tilde\alpha_k^i,\alpha_k^i\} > \min\{\bar\epsilon^p,(g_{\min})_k^2\}.
\end{equation}
By the instructions of Algorithm LOG-DFL, we have that, for all $k$ sufficiently large,
\begin{equation}\label{monotona_non_cresc}
    P(x_{k+1};\bar\eps) \leq P(y_k^{n+1};\bar\eps)\leq\dots\leq P(y_{k}^1;\bar\eps) = P(x_k;\bar\eps).
\end{equation}
Hence,
\begin{equation}\label{P_to_barP}
    \lim_{k\to\infty} P(x_k;\bar\eps) = \bar P<+\infty.
\end{equation}
Then, recalling Proposition \ref{eps_to_0_1}, $\forall k$ sufficiently large, we have that:
\begin{equation}
\label{utile0n_new} \lim_{k\to\infty}\tilde\alpha_k^i=0\quad\quad {\it
for~}i=1,\ldots,n,
\end{equation}
\begin{equation}
\label{utile1n_new} \lim_{k\to\infty}\alpha_k^i=0\quad\quad {\it
for~}i=1,\ldots,n,
\end{equation}
%
\par 
Now, recalling  (\ref{utile3n_new}) and the fact that $\eps_k = \bar \eps$ for all $k$ sufficiently large, we have that:
\[
\lim_{k\to\infty}(g_{min})_k=0,
\]
Given the definition of $(g_{min})_k$ in the algorithm and the fact that the number of constraints $m$ and of the variables $n$ are both finite, an infinite index set $K^{\prime\prime}\subseteq\{0,1,\dots\}$ exists such that 
    \[
    (g_{min})_k = |g_{\bar\jmath}(y_k^{\bar\imath})|,
    \]
for some $\bar\jmath\in\{1,\dots,m\}$ and $\bar\imath\in\{1,\dots,n+1\}$.
%
    \[
    y_k^{\bar\imath} = x_k +\sum_{\ell = 1}^{\bar\imath -1} \alpha_k^\ell d_k^\ell.
    \]
Now, since $x_k\in X$ then, a subset of indices $K^{\prime\prime\prime}\subseteq K^{\prime\prime}$ exists such that
\begin{eqnarray*}
&& \lim_{k\to\infty,k\in K^{\prime\prime\prime}} x_k = \bar x \\
&& \lim_{k\to\infty,k\in K^{\prime\prime\prime}} y_k^{\bar\imath} = \bar x.
\end{eqnarray*}
Then, we have that
\[
 \lim_{k\to\infty,k\in K^{\prime\prime\prime}} |g_{\bar\jmath}(x_k)| = 0 
\]
i.e. $\bar x\in\partial \cal S$, meaning that $P(\bar x;\bar\eps) = +\infty$. This is a contradiction with (\ref{P_to_barP}) and concludes the proof.
\end{proof}

We introduce the following index set 
\begin{equation}\label{definition_of_K}
 K = \{k:\ \eps_{k+1} < \eps_k\}.
\end{equation}
Note that, by virtue of Proposition \ref{eps_to_0}, $K$ is an infinite index set.
\par\medskip
{In the next propositions we report two technical results needed to show the convergence properties of the algorithm. The first one guarantees     the convergence to zero of the sequences of the step sizes produced by the algorithm. The second one points out that, eventually, the algorithm performs suitable samplings of the merit function along all the generators of the cone of feasible directions.}

\begin{proposition}
\label{alfa_to_0}
Let $\{\tilde\alpha_k^i\}$ and $\{\alpha_k^i\}$ be the sequences produced by  Algorithm LOG-DFL. Then,
\[
\lim_{k\to\infty}\max_{i=1,\dots,n}\{\tilde\alpha_k^i,\alpha_k^i\} = 0.
\]
\end{proposition}

\begin{proof}
The proof follows from the updating rule of the algorithm
\[
\max_{i=1,2,\dots,n}\{\tilde\alpha_k^i,\alpha_k^i\} \le \min\{\epsilon_k^p,(g_{min})_k^2\},
\]
and Proposition \ref{eps_to_0}.
\end{proof}


In the next proposition we report some technical results similar to those stated in Proposition \ref{part_a_theorem_1} and that will be used in the proof of the main convergence theorem.

\begin{proposition}\label{part_a_theorem}
Let $\{x_k\}$, $\{\eps_k\}$, and $\{y_k^i\}$, $i=1,\dots,n+1$, be the sequences produced by Algorithm LOG-DFL {and let $\{x_k\}_{\tilde K}$ be a subsequence converging to the point $\bar x$. Then,  for all $d^i\in D\cap D(\bar x)$,
there exist 
scalars $\xi_k^i>0$ such that:

\smallskip
\noindent for $k\in {\tilde K}$} sufficiently large, 
\begin{eqnarray}
 \label{d0bis_old}
 && y_k^i+ \xi_k^id^i\in X\cap \stackrel{\circ}{\cal S}, \\
 %
 \label{d1bis1_old}
 && P({y_k^i}+ \xi_k^i d^i;\eps_k)\ge P({y_k^i};\eps_k)-o(\xi_k^i);
\end{eqnarray}
and, 
\begin{eqnarray}
 %
 \label{xi_ki_to_zero}
    &&  \lim_{k\to\infty{, k\in{\tilde K}}} \xi_k^i = 0\\
\label{y_ki_to_xbar}
    && \lim_{k\to\infty{, k\in{\tilde K}}} \|y_k^i - x_k\| = 0
\end{eqnarray}
 
 \end{proposition}
\begin{proof}
We recall that, by the instructions of Algorithm DFL, at every
iteration~$k$, the following set of directions is considered:
$$
D_k =\{d_k^1,-d_k^1,\dots,d_k^n,-d_k^n\}=\{\pm e^1,\dots,\pm e^n\}=D.
$$
At every iteration $k$, Algorithm DFL extracts information on the behavior of the penalty function along both $d_k^i$ and
$-d_k^i$.

In particular, along all $d_k^i$, $i=1,\dots,n$, the algorithm identifies the following circumstances:

\begin{itemize}
    \item[i)] (Step 1.4 is executed) let us define $(\alpha_k^i)^+$, $(\alpha_k^i)^-$ such that
    \begin{eqnarray*}
     && (\alpha_k^i)^+ = 0,\ \mbox{or}\\ 
     && \qquad y_k^i + (\alpha_k^i)^+ d_k^i \in X\cap {\cal S},\quad P(y_k^i + (\alpha_k^i)^+ d_k^i;\eps_k) \geq P(y_k^i;\eps_k) - \gamma ((\alpha_k^i)^+)^2
    \end{eqnarray*}
    \begin{eqnarray*}
     && (\alpha_k^i)^- = 0,\ \mbox{or}\\ 
     && \qquad y_k^i - (\alpha_k^i)^- d_k^i \in X\cap {\cal S},\quad P(y_k^i - (\alpha_k^i)^- d_k^i;\eps_k) \geq P(y_k^i;\eps_k) - \gamma ((\alpha_k^i)^-)^2
    \end{eqnarray*}
    \item[ii)] (Expansion step is executed at Step 1.3) let us define $(\alpha_k^i)^+$ and $\check\alpha_k^i$  such that
    \begin{eqnarray*}
     && (\alpha_k^i)^+ = 0,\ \mbox{or}\\ 
     && \qquad y_k^i + (\alpha_k^i)^+ d_k^i \in X\cap {\cal S},\quad P(y_k^i + (\alpha_k^i)^+ d_k^i;\eps_k) \geq P(y_k^i;\eps_k) - \gamma ((\alpha_k^i)^+)^2
    \end{eqnarray*}
    \[
     y_k^i - \check\alpha_k^i d_k^i \in X\cap {\cal S},\quad P(y_k^i - \check\alpha_k^i d_k^i;\eps_k) \geq P(y_k^i;\eps_k) - \gamma (\check\alpha_k^i)^2
    \]    
    \item[iii)] (Expansion step is executed at Step 1.2) let us define $\tilde{y}_k^i = y_k^i + \alpha_k^i d_k^i$.\\ 
    $\alpha_k^i$ and $\check\alpha_k^i$ such that
    \[
     y_k^i + \alpha_k^i d_k^i \in X\cap {\cal S},\quad P(\tilde{y}_k^i;\eps_k) \geq P(\tilde{y}_k^i - \alpha_k^id_k^i;\eps_k) - (-\gamma (\alpha_k^i)^2)
    \]
    \[
     y_k^i + \check\alpha_k^i d_k^i \in X\cap {\cal S},\quad P(y_k^i + \check\alpha_k^i d_k^i;\eps_k) \geq P(y_k^i;\eps_k) - \gamma (\check\alpha_k^i)^2
    \]
\end{itemize}

Furthermore, recalling Proposition \ref{alfa_to_0} , we also have that
\begin{equation}
\label{utile1tris}
\lim_{k\to\infty}\alpha_k^i=0\quad\quad {\rm
for~}i=1,\ldots,n.
\end{equation}
\begin{equation}
\label{utile1bisntris}
\lim_{k\to\infty}\tilde\alpha_k^i=0\quad\quad {\rm
for~}i=1,\ldots,n.
\end{equation}
Then, since $(\alpha_k^i)^+ \leq \tilde\alpha_k^i$, $(\alpha_k^i)^- \leq \tilde\alpha_k^i$, $\check\alpha_k^i \leq \frac{\alpha_k^i}{\delta}$, we also have that
\begin{equation}\label{alpha+tozero}
\lim_{k\to\infty}(\alpha_k^i)^+=0\quad\quad {\rm
for~}i=1,\ldots,n.
\end{equation}
\begin{equation}\label{alpha-tozero}
\lim_{k\to\infty}(\alpha_k^i)^-=0\quad\quad {\rm
for~}i=1,\ldots,n.
\end{equation}
\begin{equation}\label{checkalphatozero}
\lim_{k\to\infty}\check\alpha_k^i=0\quad\quad {\rm
for~}i=1,\ldots,n.
\end{equation}
By recalling the definitions of the search direction $d_k^i$,
$i=1,\ldots,n$, we obtain
\begin{equation}
\label{inclusione1}
D\cap D(\bar x)\subseteq \{ d_k^1,-d_k^1,\dots,d_k^n,-d_k^n\}.
\end{equation}
Now by using (\ref{utile1bisntris}), (\ref{inclusione1}) and Proposition \ref{iota}, we have that, for sufficiently large $k \in \tilde K$ and for all $d_k^i \in D\cap D(\bar x)$, $(\alpha_k^i)^+=0$ can not happen and that, for sufficiently large $k \in \tilde K$ and for all $-d_k^i \in D\cap D(\bar x)$, $(\alpha_k^i)^-=0$ can not happen.

Let us consider all the directions
$d^i\in D\cap D(\bar x)$.

If $d^i=d_k^i$, by setting $\xi_k^i=(\alpha_k^i)^+$, and
$o(\xi_k^i)=\gamma ((\alpha_k^i)^+)^2$ if we are in i) or ii); otherwise by setting $\xi_k^i=\check\alpha_k^i$, and
$o(\xi_k^i)=\gamma (\check\alpha_k^i)^2$\ if 
we are in iii), for sufficiently large $k\in \tilde K$, we can write
\begin{equation}
\label{d0bis}
y_k^i+ \xi_k^id^i\in X\cap \stackrel{\circ}{\cal S},
\end{equation}
\begin{equation}
\label{d1bis1}
P({y_k^i}+ \xi_k^i d^i;\eps_k)\ge P({y_k^i};\eps_k)-o(\xi_k^i),
\end{equation}
On the other hand, if $d^i=-d_k^i$, (\ref{d0bis}), (\ref{d1bis1}) hold, for sufficiently large $k\in \tilde K$, by
setting $\xi_k^i=(\alpha_k^i)^-$, and
$o(\xi_k^i)=\gamma ((\alpha_k^i)^-)^2$ if we are in case i);
by setting $\xi_k^i=\check\alpha_k^i$ and
$o(\xi_k^i)=\gamma (\check\alpha_k^i)^2$ if we are in case ii); by setting $\xi_k^i=\alpha_k^i$, $y_k^i=\tilde{y}_k^i$ and
$o(\xi_k^i)=-\gamma (\alpha_k^i)^2$ if we are in case iii).\par
Then, given the definition of the scalars $\xi_k^i$, we have that (\ref{xi_ki_to_zero}) is satisfied.\par
Finally, since $y_k^i = x_k + \sum_{j=1}^{i-1}\alpha_k^j d_k^j$, recalling (\ref{utile1tris}), we obtain that (\ref{y_ki_to_xbar}) is also satisfied.
\end{proof}

Finally it is possible to state the main result concerning the convergence properties of the proposed algorithm.

\begin{theorem}
\label{mainalgdfl} Let $\{x_k\}$ be the sequence generated by
Algorithm {\rm DFL}. Let $K$ be the set of indices defined in (\ref{definition_of_K}). Assume that every limit point of the sequence
$\{x_k\}_K$ satisfies the MFCQ; then, every limit point
$\bar x$ of the subsequence $\{x_k\}_K$ is a stationary point of
problem {\rm (\ref{prob1})}.
\end{theorem}

\begin{proof}
Since $\{x_k\}_K\subseteq X$ and $X$ is compact, the subsequence $\{x_k\}_K$ admits limit points. Let us consider one such limit point $\bar x$, i.e. an index set $\bar K\subseteq K$ exists such that 
\[
\lim_{k\to\infty,k\in \bar K} x_k = \bar x.
\]
Let us denote $\bar D=D\cap D(\bar x)$. {Recalling Proposition \ref{part_a_theorem}
we have that (\ref{d0bis_old}), (\ref{d1bis1_old}),   (\ref{xi_ki_to_zero}) and (\ref{y_ki_to_xbar}) hold.}

\par\medskip\noindent
We have that for all $k\in\bar K$ sufficiently large and for all $d^i\in D\cap{D}(\bar{x})$:
\begin{itemize}
    \item[i)] $y_k^i\in X\cap \stackrel{\circ}{\cal S}$
    \item[ii)]$y_k^i+\xi_k^id^i\in X\cap \stackrel{\circ}{\cal S}$
\end{itemize}
Now, let us assume there exists $\hat{t}_k^i \in (0,1)$ such that $y_k^i+\hat{t}_k^i \xi_k^i d^i \notin  X\cap \stackrel{\circ}{\cal S}$. Then, by the compactness of $X$
\[
y_k^i+\hat{t}_k^i\xi_k^id^i \notin \stackrel{\circ}{\cal S}, \quad \text{i.e. } g_l(y_k^i+\hat{t}_k^i\xi_k^id^i)>0
\]
Using the continuity assumption on the constraints, there exists at least one constant $\check{t}_k^i \in (0,\hat{t}_k^i)$ such that:
\[
g_l(y_k^i+\check{t}_k^i \xi_k^id^i)=0,
\]
if multiple constants exist that satisfy the condition above, we will consider $\check{t}_k^i$ to be the smallest one. We have now that $y_k^i + t\xi_k^id^i \in X\cap \stackrel{\circ}{\cal S}$ for all $t \in [0,\check{t}_k^i)$, so that by the definition of $P$:
\[
P(y;\eps_k) \text{ is continuous } \forall y \in [y_k^i,y_k^i+\check{t}_k^i\xi_k^id^i) \text{ and } P(y_k^i+\check{t}_k^i\xi_k^id^i;\eps_k)=+\infty,
\]
thus, a constant $t_k^{*,i} \in (0,\bar{t}_2)$ must exist such that:
\begin{eqnarray}
    && \label{new_fallimento}P({y_k^i}+ t_k^{*,i} \xi_k^i d^i;\eps_k)\ge P({y_k^i};\eps_k)-o(t_k^{*,i}\xi_k^i)\\
    && \label{new_feasibility}y_k^i+t \xi_k^i d^i \in X\cap \stackrel{\circ}{\cal S} \forall t\in[0,t_k^{*,i}].
\end{eqnarray}
Let us now set
\[
\bar{\xi_k^i} = t_k^{*,i} \xi_k^i \le \xi_k^i
\]
and apply the mean-value theorem to (\ref{new_fallimento}). Hence, we can write
$$
-o(\bar\xi_k^{i})\le P(y_k^i+\bar\xi_k^{i} d^{i};\eps_k) -
P(y_k^i;\eps_k)= \bar\xi_k^{i}\nabla
P(u_k^i;\eps_k)^{T}d^{i}\quad\quad \forall d^i\in \bar D,
$$
where $u_k^i=y_k^i+t_k^{i}\bar\xi_k^{i} d^{i}$, with $t_k^{i}\in (0,1)$. By recalling (\ref{new_feasibility}), $u_k^i \in int({\cal S})$. Thus, we have
\[
-\frac{o(\bar\xi_k^{i})}{\bar\xi_k^{i}} \leq \nabla P(u_k^i;\eps_k)^{T}
d^i\quad\quad \forall d^i\in \bar D.
\]
By considering the expression of $P(x;\eps)$, we can write
\begin{eqnarray}
\label{prop2.eq1old}
&&{\nabla P(u_k^i;\eps_k)^T d^i=}  \Bigg(\nabla
f(u_k^i)+\sum_{l=1}^m\frac{\eps_k}{-g_l(x)}\nabla
g_l(u_k^i) \\\nonumber 
&& \qquad\qquad\qquad+ \sum_{j=1}^q \frac{\nu}{\eps_k} |h_j(u_k^i)|^{\nu-1}\nabla h_j(u_k^i) \Bigg)^{T} d^{i} \geq
-\frac{o(\bar\xi_k^{i})}{\bar\xi_k^{i}}\quad\quad \forall d^i\in \bar D.
\end{eqnarray}
Recalling that $u_k^i=y_k^i+t_k^{i}\bar\xi_k^{i} d^{i}$, with
$t_k^{i}\in (0,1)$, and that $\bar\xi_k^i\leq\xi_k^i$ for all $i$, we have that, for all $i$ such that $d^i\in
\bar D$,
\begin{equation}
\label{prop2.eq5sec} \lim_{k\rightarrow\infty, k\in \bar K} u_k^i =
\bar x.
\end{equation}

\par\medskip
\par\medskip
{Now it is possible to define the following approximations of the multipliers.

For $l=1,\dots,m$ set
$$
\la_l(x;\eps) = \frac{\eps}{-g_l(x)}
$$

For $j=1,\dots,q$ set
$$
\mu_j(x;\eps) = \frac{\nu}{\eps} \ |h_j(x)|^{\nu - 1}
$$
The sequences $\{\la_l(x_k;\eps_k)\}_{\bar K}$, $l=1,\dots,m$, and $\{\mu_j(x_k;\eps_k)\}_{\bar K}$, $j=1,\dots,q$ are bounded. The proof of this property  is rather technical and, to simplify the exposition, it is reported in the appendix as Proposition~\ref{theo3-lambda}.
\par\medskip
Then there exists a subset of $\bar K$,
which we relabel again $\bar K$, such that
\[
\lim_{k\to\infty, k\in \bar K}\la_l(x_k;\eps_k) = \bar\la_l\geq
0,\qquad l=1,\dots,m,
\]
\[
\lim_{k\to\infty, k\in \bar K}\mu_j(x_k;\eps_k) = \bar\mu_j\geq
0,\qquad j=1,\dots,q,
\]
where $\bar\la_l = 0$ for $l\not\in I(\bar x)$.}
\par\medskip

Since $y_k^i \in \cal S$ and, by continuity of $g_\ell(x)$ for all $\ell=1,\dots,m$, $\cal S$ is closed, any accumulation point of $\{y_k^i\} \in \cal S$.
We consider now the sequence of positive penalty
parameters ${\epsilon_k}$.
By~ Proposition \ref{eps_to_0}, we have that:
$$
\lim_{k\to\infty}\epsilon_k=0
$$
recalling assumption (i), recalling the continuity assumptions, multiplying ~(\ref{prop2.eq1old}) by $\epsilon_k$ and taking the limit, we have:
\begin{equation}
\left(\sum_{j=1}^p \nu |h_j(\xbar)|^{\nu -1}\nabla h_j(\xbar)\right)^{T} d^{i} \geq 0 \ \ \ \ \ \forall d^i\in \bar D.
\end{equation}
Since $\xbar$ satisfies MFCQ, by (\ref{box_mfcq}), it must result:
\[
h_j(\xbar)=0 \qquad \forall  j=1,\dots,p.
\]
Therefore the point $\xbar$ is feasible.
\par\medskip
 By simple manipulations, (\ref{prop2.eq1old}) can be rewritten as
\begin{eqnarray}
\label{prop2.eq3-1old}
\lefteqn{\Bigg( \nabla f(u_k^i) + \sum_{l=1}^m\nabla g_l(u_k^i)\la_l(x_k;\eps_k) } \\\nonumber
 & + & \sum_{l=1}^m\nabla
g_l(u_k^i)\left(\la_l(u_k^i;\eps_k)-\la_l(x_k;\eps_k)\right)+\sum_{j=1}^q \nabla h_j(u_k^i)\mu_j(x_k;\eps_k) \\\nonumber
& + & \sum_{j=1}^1\nabla
h_j(u_k^i)\left(\mu_j(u_k^i;\eps_k)-\mu_j(x_k;\eps_k)\right) \Bigg)^T d^{i} \geq -\frac{o(\xi_k^{i})}{\xi_k^{i}}\quad\quad
\forall i: d^i\in \bar D.
\end{eqnarray}
Taking the limits for $k\to\infty$ and $k\in \bar K$ in relation (\ref{prop2.eq3-1old}) and recalling (\ref{limite_lambda}) and (\ref{limite_mu}) from the proof of
Proposition~\ref{theo3-lambda} previously invoked, we obtain
\[
\left(\nabla f(\bar x) + \sum_{l=1}^m\nabla g_l(\bar x)\bar\la_l + \sum_{j=1}^q\nabla h_j(\bar x)\bar\mu_j
\right)^T d^i \geq 0\qquad \forall i:d^i\in \bar D.
\]
Recalling that $\bar D=D\cap D(\bar x)$, from Proposition~\ref{cn1} we get
\[
\nabla L(\bar x,\bar\la,\bar\mu)^T d \geq 0\qquad \forall d\in D(\bar x),
\]
which concludes the proof.
\end{proof}

\section{Numerical experiments}
\label{sec:experiments}

In this section we report the numerical performance of the
proposed log-barrier derivative-free Algorithm LOG-DFL on
a set of test problems chosen from a well-known collection.

\subsection{Test problem collection}
In this subsection we report the set of constrained test problems selected from the CUTEst collection \cite{gould2015cutest}. In particular, we selected all the problems with $n\leq 50$ variables and having at least one inequality constraint for which the provided initial point is strictly feasible, i.e. such that at least an index $j$ exists with $g_j(x_0) < 0$ {(the constraints such that $g_j(x_0) \ge  0$ are taken into account by an exterior penalty term)}. This gives us a total of $96$ problems.

\begin{table}[ht]
\begin{center}\scriptsize
\begin{minipage}{0.48\textwidth}
\begin{center}
\begin{tabular}{|l|c|c|c|}\hline
Problem & $n_p$ & $m_p$ & $\bar m_p$ \\\hline
   ANTWERP & 27 &  10 &   2 \\
    DEMBO7 & 16 &  21 &  16 \\
  ERRINBAR & 18 &   9 &   1 \\
     HS117 & 15 &   5 &   5 \\
     HS118 & 15 &  29 &  28 \\
    LAUNCH & 25 &  29 &  20 \\
   LOADBAL & 31 &  31 &  20 \\
   MAKELA4 & 21 &  40 &  20 \\
      MESH & 33 &  48 &  17 \\
  OPTPRLOC & 30 &  30 &  28 \\
       RES & 20 &  14 &   2 \\
  SYNTHES2 & 11 &  15 &   1 \\
  SYNTHES3 & 17 &  23 &   1 \\
  TENBARS1 & 18 &   9 &   1 \\
  TENBARS4 & 18 &   9 &   1 \\
  TRUSPYR1 & 11 &   4 &   1 \\
  TRUSPYR2 & 11 &  11 &   8 \\
      HS12 &  2 &   1 &   1 \\
      HS13 &  2 &   1 &   1 \\
      HS16 &  2 &   2 &   2 \\
      HS19 &  2 &   2 &   1 \\
      HS20 &  2 &   3 &   3 \\
      HS21 &  2 &   1 &   1 \\
      HS23 &  2 &   5 &   4 \\
      HS30 &  3 &   1 &   1 \\
      HS43 &  4 &   3 &   3 \\
      HS65 &  3 &   1 &   1 \\
      HS74 &  4 &   5 &   2 \\
      HS75 &  4 &   5 &   2 \\
      HS83 &  5 &   6 &   5 \\
      HS95 &  6 &   4 &   3 \\
      HS96 &  6 &   4 &   3 \\
      HS97 &  6 &   4 &   2 \\
      HS98 &  6 &   4 &   2 \\
     HS100 &  7 &   4 &   4 \\
     HS101 &  7 &   6 &   2 \\
     HS104 &  8 &   6 &   3 \\
     HS105 &  8 &   1 &   1 \\
     HS113 & 10 &   8 &   8 \\
     HS114 & 10 &  11 &   8 \\
     HS116 & 13 &  15 &  10 \\
      S365 &  7 &   5 &   2 \\
   ALLINQP & 24 &  18 &   9 \\
  BLOCKQP1 & 35 &  16 &   1 \\
  BLOCKQP2 & 35 &  16 &   1 \\
  BLOCKQP3 & 35 &  16 &   1 \\
  BLOCKQP4 & 35 &  16 &   1 \\
  BLOCKQP5 & 35 &  16 &   1 \\\hline
\end{tabular}
\end{center}
\end{minipage}
\begin{minipage}{0.48\textwidth}
\begin{center}
\begin{tabular}{|l|c|c|c|}\hline
Problem & $n_p$ & $m_p$ & $\bar m_p$ \\\hline
  CAMSHAPE & 30 &  94 &  90 \\
      CAR2 & 21 &  21 &   5 \\
  CHARDIS1 & 28 &  14 &  13 \\
       EG3 & 31 &  90 &  60 \\
  GAUSSELM & 29 &  36 &  11 \\
       GPP & 30 &  58 &  58 \\
  HADAMARD & 37 &  93 &  36 \\
   HANGING & 15 &  12 &   8 \\
  JANNSON3 & 30 &   3 &   2 \\
  JANNSON4 & 30 &   2 &   2 \\
   KISSING & 37 &  78 &  32 \\
  KISSING1 & 33 & 144 & 113 \\
  KISSING2 & 33 & 144 & 113 \\
  LIPPERT1 & 41 &  80 &  64 \\
  LIPPERT2 & 41 &  80 &  64 \\
   LUKVLI1 & 30 &  28 &  28 \\
  LUKVLI10 & 30 &  28 &  14 \\
  LUKVLI11 & 30 &  18 &   3 \\
  LUKVLI12 & 30 &  21 &   6 \\
  LUKVLI13 & 30 &  18 &   3 \\
  LUKVLI14 & 30 &  18 &  18 \\
  LUKVLI15 & 30 &  21 &   7 \\
  LUKVLI16 & 30 &  21 &  13 \\
  LUKVLI17 & 30 &  21 &  21 \\
  LUKVLI18 & 30 &  21 &  21 \\
   LUKVLI2 & 30 &  14 &   7 \\
   LUKVLI3 & 30 &   2 &   2 \\
   LUKVLI4 & 30 &  14 &   4 \\
   LUKVLI6 & 31 &  15 &  15 \\
   LUKVLI8 & 30 &  28 &  14 \\
   LUKVLI9 & 30 &   6 &   6 \\
     MANNE & 29 &  20 &  10 \\
  MOSARQP1 & 36 &  10 &  10 \\
  MOSARQP2 & 36 &  10 &  10 \\
     NGONE & 29 & 134 & 106 \\
  NUFFIELD & 38 & 138 &  28 \\
   OPTMASS & 36 &  30 &   6 \\
   POLYGON & 28 & 119 &  94 \\
  POWELL20 & 30 &  30 &  15 \\
  READING4 & 30 &  60 &  30 \\
  SINROSNB & 30 &  58 &  29 \\
  SVANBERG & 30 &  30 &  30 \\
  VANDERM1 & 30 &  59 &  29 \\
  VANDERM2 & 30 &  59 &  29 \\
  VANDERM3 & 30 &  59 &  29 \\
  VANDERM4 & 30 &  59 &  29 \\
       YAO & 30 &  30 &   1 \\
    ZIGZAG & 28 &  30 &   5 \\\hline
\end{tabular}
\end{center}
\end{minipage}
\end{center}
\caption{Set of test problems selected from the CUTEst collection. $n_p$, $m_p$, and $\bar m_p$ denote, respectively, the number of variables, of constraints and of strictly feasible inequality constraints for the given problem.}
\end{table}

\begin{figure}
    \centering
    \includegraphics[width=0.48\textwidth]{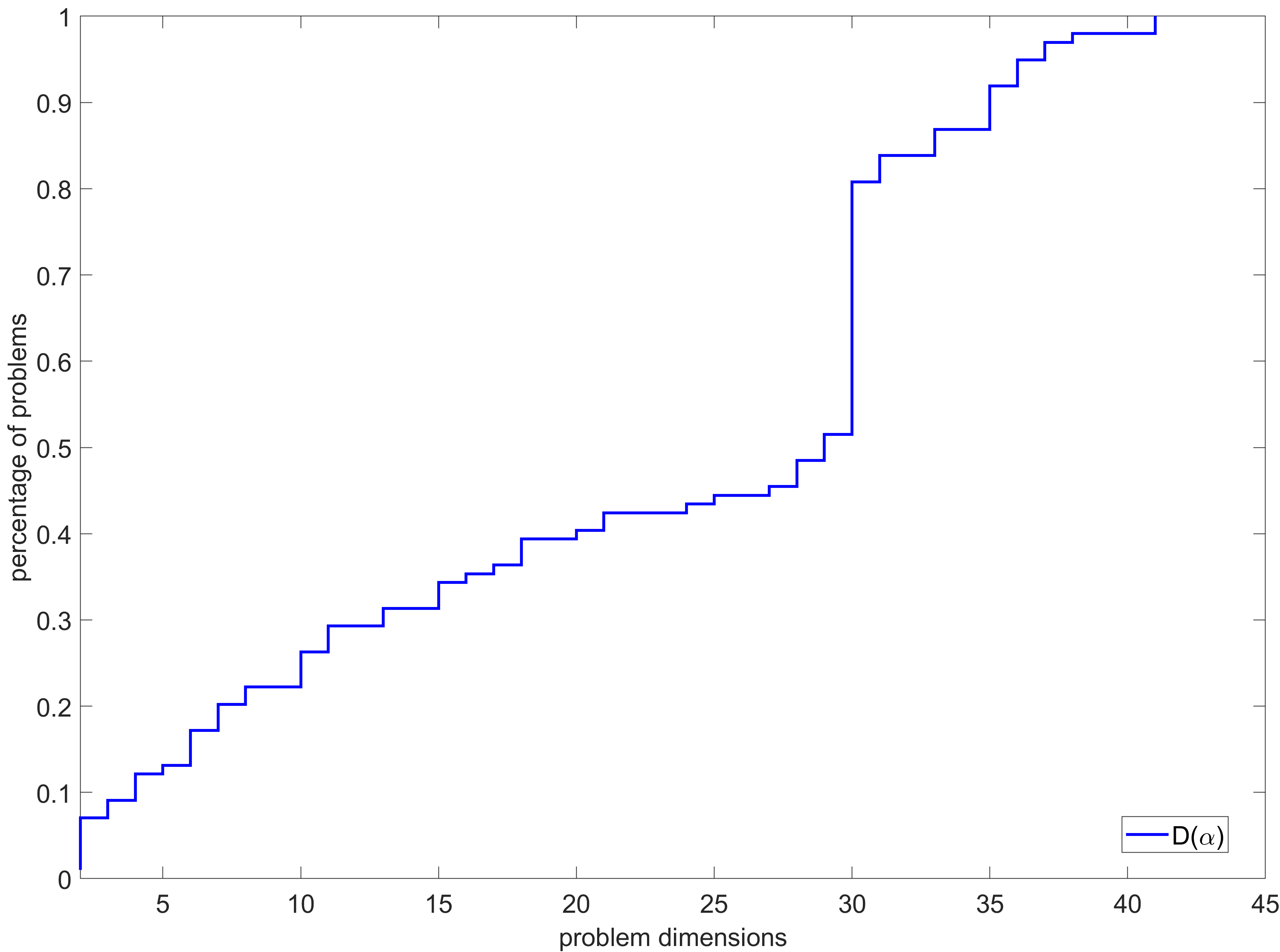}
    \includegraphics[width=0.48\textwidth]{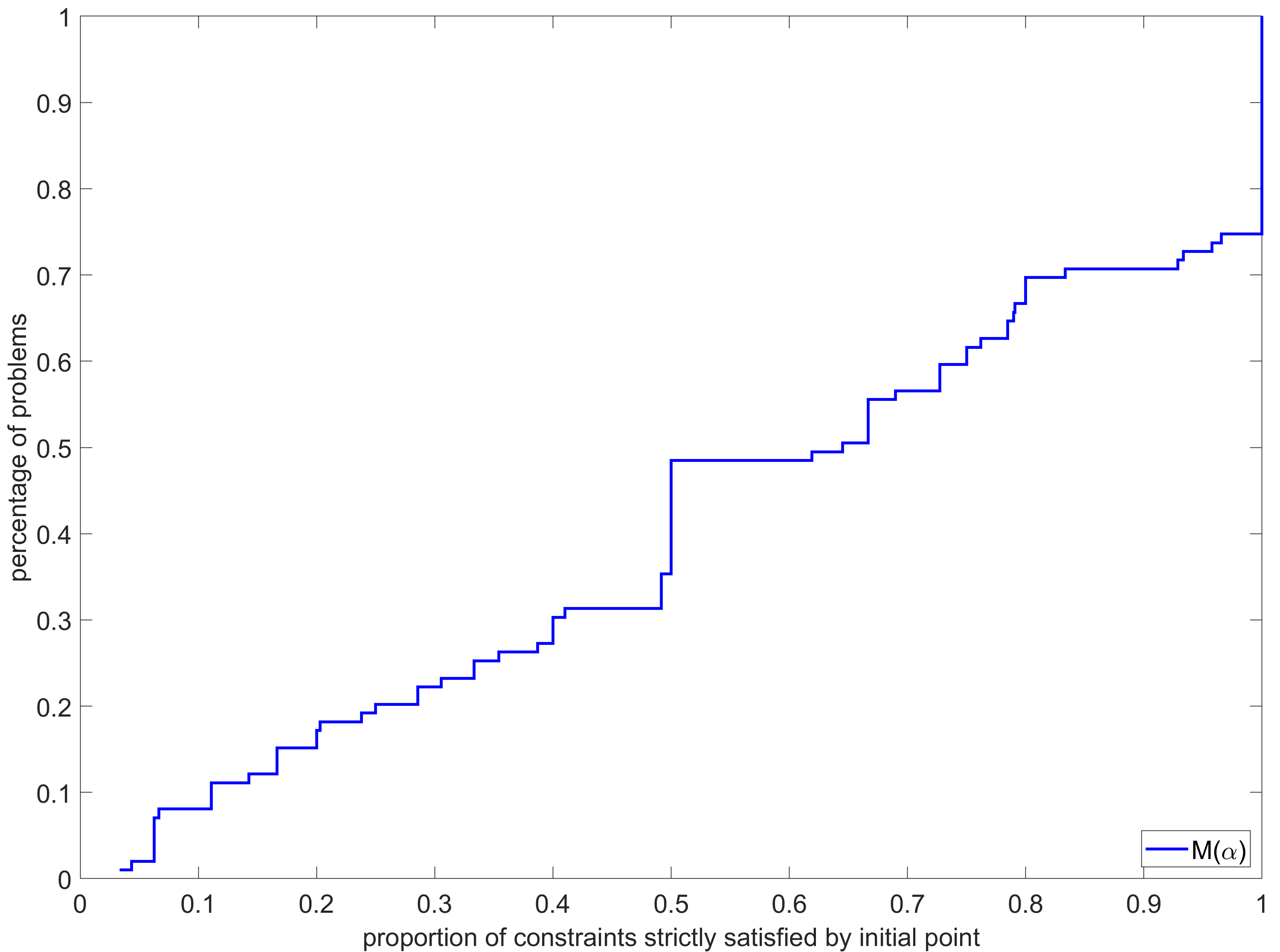}
    \caption{Cumulative distributions, respectively, of the number of variables and of the proportion of strictly satisfied constraints with respect to the total number of constraints.}
    \label{fig:problems}
\end{figure}

In figure \ref{fig:problems} we report the cumulative distributions, respectively, of the number of variables and of the proportion of strictly satisfied constraints with respect to the total number of constraints, i.e.
\begin{eqnarray*}
D(\alpha) & = & \frac{1}{N}\left|\left\{p\in{\cal P}:\ n_p\leq \alpha\right\}\right| \\
M(\alpha) & = & \frac{1}{N}\left|\left\{p\in{\cal P}:\ \frac{\bar m_p}{m_p}\leq \alpha\right\}\right|
\end{eqnarray*}
where
\begin{itemize}
    \item[-] ${\cal P}$ is the set of problems;
    \item[-] $N = |{\cal P}|$;
    \item[-] $n_p$ is the number of variables of problem $p\in\cal P$;
    \item[-] $m_p$ is the number of constraints of problem $p\in\cal P$;
    \item[-] $\bar m_p$ is the number of strictly satisfied inequality constraints at the initial point for problem $p\in\cal P$.
\end{itemize}








\subsection{Implementation details}

The proposed method has been implemented in Python, and all the experiments have been conducted by
choosing the following:
\begin{itemize}
    \item[-] {\em Exponent parameter for exterior penalty}\par $\nu=1.1$. \hfill \break
    \item[-] {\em Parameters introduced in the LOG-DFL algorithm  scheme}\par $\gamma = 10^{-4}$, $p=1.1$, $\tilde\alpha^i_0 = 1.0$. As concerns the parameter $\theta$, we split it into two different parameters: $\theta^{in}=0.35$ and $\theta^{ex}=10^{-2}$. Their role will be explained in the following {\em penalty parameter updating criterion} point.  \hfill \break
    \item[-] {\em Penalty parameters initialization}\par
    We computed the values of the constraints at the starting point $x_0$ and we defined two sets of indices:
    \[
    J_{in} := \{j: g_j(x_0)<0\},
    \]
    \[
    J_{ex} := \{j: g_j(x_0) \ge 0\}.
    \]
    So we set a log penalty for feasible constraints and exterior penalty for active and unfeasible ones. We define two parameters $\eps_0^{in}$, $\eps_0^{ex}$, and we initialize it with the following rules:
    \[
    \eps_0^{in} = 10^{-1}
    \]
    \[
    \eps_0^{ex} = \min\left\{10^{-1}, \frac{1}{|f(x_0)|}\right\}
    \]
    We can write the penalized function:
\begin{eqnarray}
\label{hybridpen}
\lefteqn{P(x,\eps_k) = } \\\nonumber
& f(x) &\ - \epsilon_k^{in} \displaystyle{\sum_{j\in J_{in}}} \log\left[-g_j(x)\right] + \frac{1}{\epsilon_k^{ex}}\left({\sum_{j \in J_{ex}}}[g_j^+(x)]^\nu + {\sum_{j=1}^q} |h_j(x)|^\nu\right),
\end{eqnarray}
    where $g_j^+(x) = \max\{g_j(x), 0\}$. \hfill \break
    \item[-] {\em Penalty parameter updating criterion}\par
    As one can see in (\ref{hybridpen}), We are now using a hybrid approach, where some inequality constraints, which will be considered as the non-relaxable ones, are handled by interior penalty and some are handled by exterior penalty. In fact, equality constraints are treated by splitting them into two opposite inequality constraints, which will be assigned with exterior penalty. According to our theoretical results and those in \cite{liuzzi2010sequential}, we use two different updating criteria:
    \begin{eqnarray}
    \label{exterior_criterion}
    \max_{i=1,2,\dots,n}\{\tilde\alpha_k^i,\alpha_k^i\} \leq (\eps_k^{ex})^p,
    \end{eqnarray}
    \begin{eqnarray}
    \label{interior_criterion}
    \max_{i=1,2,\dots,n}\{\tilde\alpha_k^i,\alpha_k^i\} \leq\min\left\{(\eps_k^{in})^p,(g_{\min})^2_k\right\},
    \end{eqnarray}
    where we remind that $(g_{\min})_k$ is the minimum absolute value for the non-relaxable constraints at iteration $k$.\\
    When (\ref{exterior_criterion}) is satisfied, the algorithm performs the following update:
    \[
    \eps_{k+1}^{ex} = \theta^{ex}\eps_k^{ex},
    \]
    When (\ref{interior_criterion}) is satisfied, the algorithm performs the following update:
    \[
    \eps_{k+1}^{in} = \theta^{in}\eps_k^{in},
    \]\hfill \break
    \item[-] {\em Stopping criterion}\par
    We stop the algorithm whenever $\max_{i=1,2,\dots,n}\{\tilde\alpha_k^i,\alpha_k^i\} \leq 10^{-14}$. Finally, we allow a maximum of 20000 function evaluations. \hfill \break
    
\end{itemize}
The LOG-DFL algorithm is freely available for download through the DFL library as package {\tt LOGDFL} at the URL \href{http://www.iasi.cnr.it/~liuzzi/DFL/}{http://www.iasi.cnr.it/$\sim$liuzzi/DFL/} 

\noindent
For comparison we used the state-of-the-art MADS algorithm implemented in  
the well-known NOMAD package (version 3.9.1) \cite{NOMAD}. NOMAD has been run using its default settings except for the type of constraints. Indeed, we forced NOMAD to handle constraints that are strictly satisfied at the initial point with an extreme barrier approach. For all other constraints, we use the mixed progressive-extreme barrier approach which is the type of management which is suggested by the developpers of NOMAD itself. Hence, for constraint $g_j$, $j=1,\dots,m$, we specified\par\medskip
\begin{center}
\texttt{EB} if $g_j(x_0) < 0$, \quad \texttt{PEB} otherwise.
\end{center}

\subsection{Performance and data profiles}
\label{sec5.1}
Results are reported in terms of performance \cite{dolan2002benchmarking} and data \cite{wild:2009} profiles which are briefly recalled in the following. In particular, let $\cal S$ be the set of solvers to be compared against each other. For each $s\in S$ and $p \in P$, the number of function evaluations required by algorithm $s$ to satisfy the convergence condition on problem $p$ is denoted as $t_{p,s}$. Given a tolerance $0< \tau < 1$ and denoted as $f_L$ the smallest objective function value computed by any algorithm on problem $p$ within a given number of function evaluations, the convergence test is 
\[
f(x_k) \leq f_L + \tau(\hat f(x_0) - f_L),
\]
where $\hat f(x_0)$ is the objective function value of the worst feasible point determined by all
the solvers (note that in the bound-constrained case, $\hat f(x_0) = f (x_0)$). The above convergence test 
requires the best point to achieve a sufficient reduction with respect to the value $\hat f(x_0)$ of the objective function at the starting point. We set to $+\infty$ the value of the objective function at infeasible points, i.e. points that have a feasibility violation $c(x) > 10^{-4}$, where 
\[
c(x) = \sum_{i=1}^m\max\{0,g_i(x)\} + \sum_{j=1}^q|h_j(x)|.
\]
Note that the smaller the value of the tolerance $\tau$ is, the higher accuracy is required at the best point. In particular, three levels of accuracy are considered in this paper for the parameter $\tau$, namely, $\tau \in \{10^{-1}, 10^{-3},10^{-5}\}$.

Performance and data profiles of solver $s$ can be formally defined as follows
\begin{eqnarray*}
	\rho_s(\alpha) & = & \frac{1}{|P|}\left|\left\{p\in P: \frac{t_{p,s}}{\min\{t_{p,s'}:s'\in S\}}\leq\alpha\right\}\right|,\\
	d_s(\kappa) & = & \frac{1}{|P|}\left|\left\{p\in P: t_{p,s}\leq\kappa(n_p+1)\right\}\right|,
\end{eqnarray*}
where $n_p$ is the dimension of problem $p$. While $\alpha$ indicates that the number of function evaluations required by algorithm $s$ to achieve the best solution is $\alpha$--times the number of function evaluations needed by the best algorithm, 
$\kappa$ denotes the number of simplex gradient estimates, with $n_p + 1$ being the number of function evaluation required to obtain one simplex gradient. Important features for the comparison are $\rho_s(1)$, which is a measure of the efficiency of the algorithm, since it is the percentage of problems for which the algorithm $s$ performs the best, and the height reached by each profile as the value of $\alpha$ or $\kappa$ increases, which measures the reliability (or robustness) of the solver, i.e. the percentage of the problems that the given solver is able to solve no matter how efficiently.

\subsection{Results and comparison}

%
%
%
%

In the following, we first describe a couple of heuristics that helped us improve our LOG-DFL algorithm.
\begin{enumerate}
\item First, drawing inspiration from the literature on interior point methods \cite{wright1997primal,forsgren2002interior,polik2010interior}, we modified our LOG-DFL algorithm by adding a mechanism that exploits a further direction $d_\rho$ which is defined using two consecutive points where the interior penalty parameter has been updated. Since the interior penalty parameter is updated when a ``quasi'' stationary point of the merit function is obtained, there are good chances that direction $d_\rho$ is a good descent direction. Moreover, $d_\rho$ should point in the direction of the ``central path'' more or less followed by the algorithm. Along such a direction we perform a further exploration by means of a suitable modification of the expansion step procedure. Such a modification is needed to satisfy (through projections) the bound constraints on the variables.  
\item Second, we try to mimic the behavior of the progressive-extreme handling of constraints ({\tt PEB} constraint type) in NOMAD. Particularly, when an initially infeasible constraint (which is penalized by using the exterior penalty approach) becomes feasible, we change the penalization method adopted for that particular constraint, thus switching to an interior point penalization. 
\end{enumerate}
LOG-DFL algorithm with heuristics is freely available for download through the DFL library as package {\tt LOGDFL} at the URL \href{http://www.iasi.cnr.it/~liuzzi/DFL/}{http://www.iasi.cnr.it/$\sim$liuzzi/DFL/} 

\begin{figure}[!ht]
    \centering
    \begin{subfigure}[b]{0.9\textwidth}
        \centering
        \includegraphics[width=0.9\textwidth]{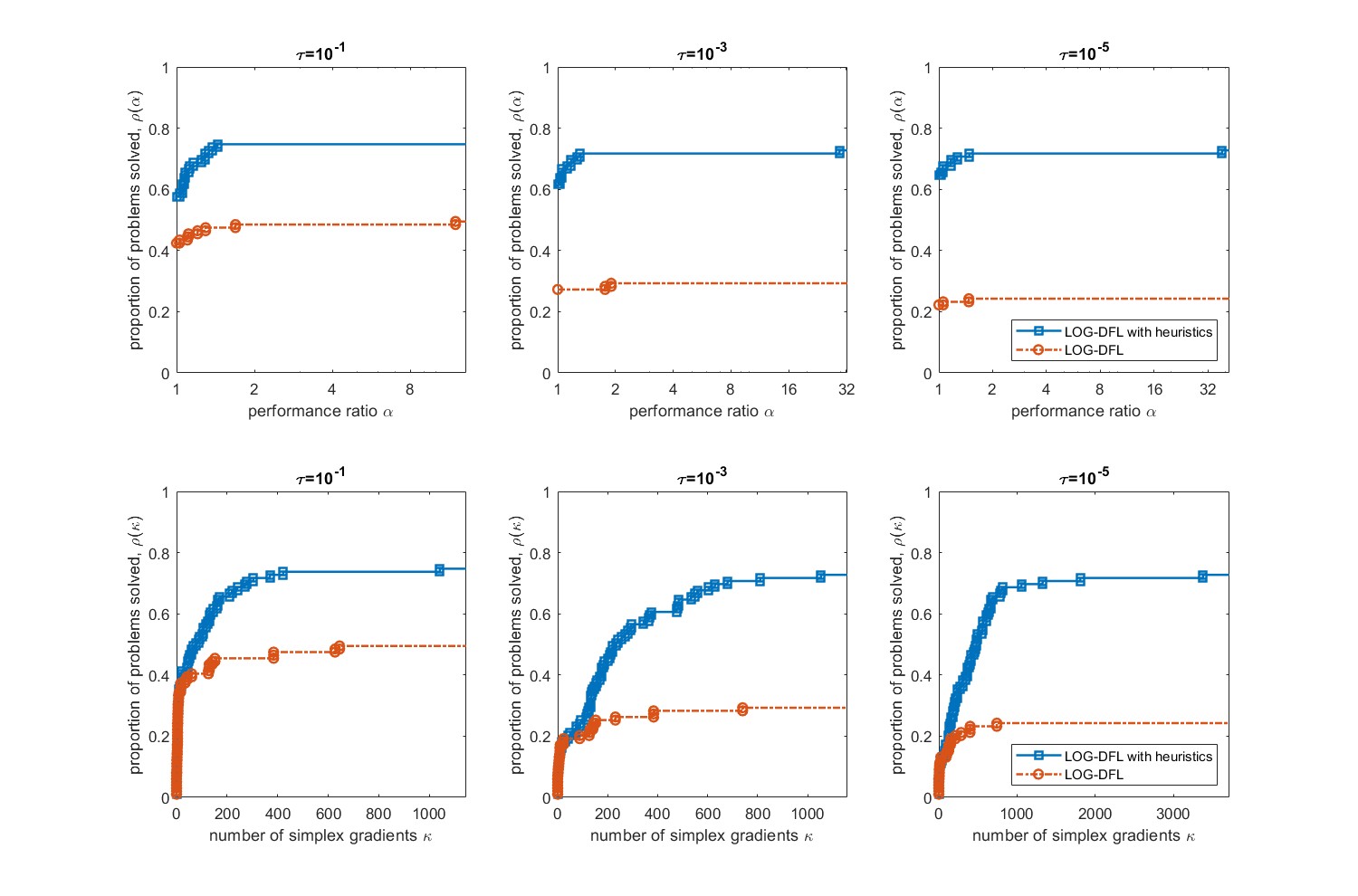}
        \caption{Whole test problems collection}
        \label{fig:restutti_innercompare}
    \end{subfigure}
    \begin{subfigure}[b]{0.9\textwidth}
        \centering
        \includegraphics[width=0.9\textwidth]{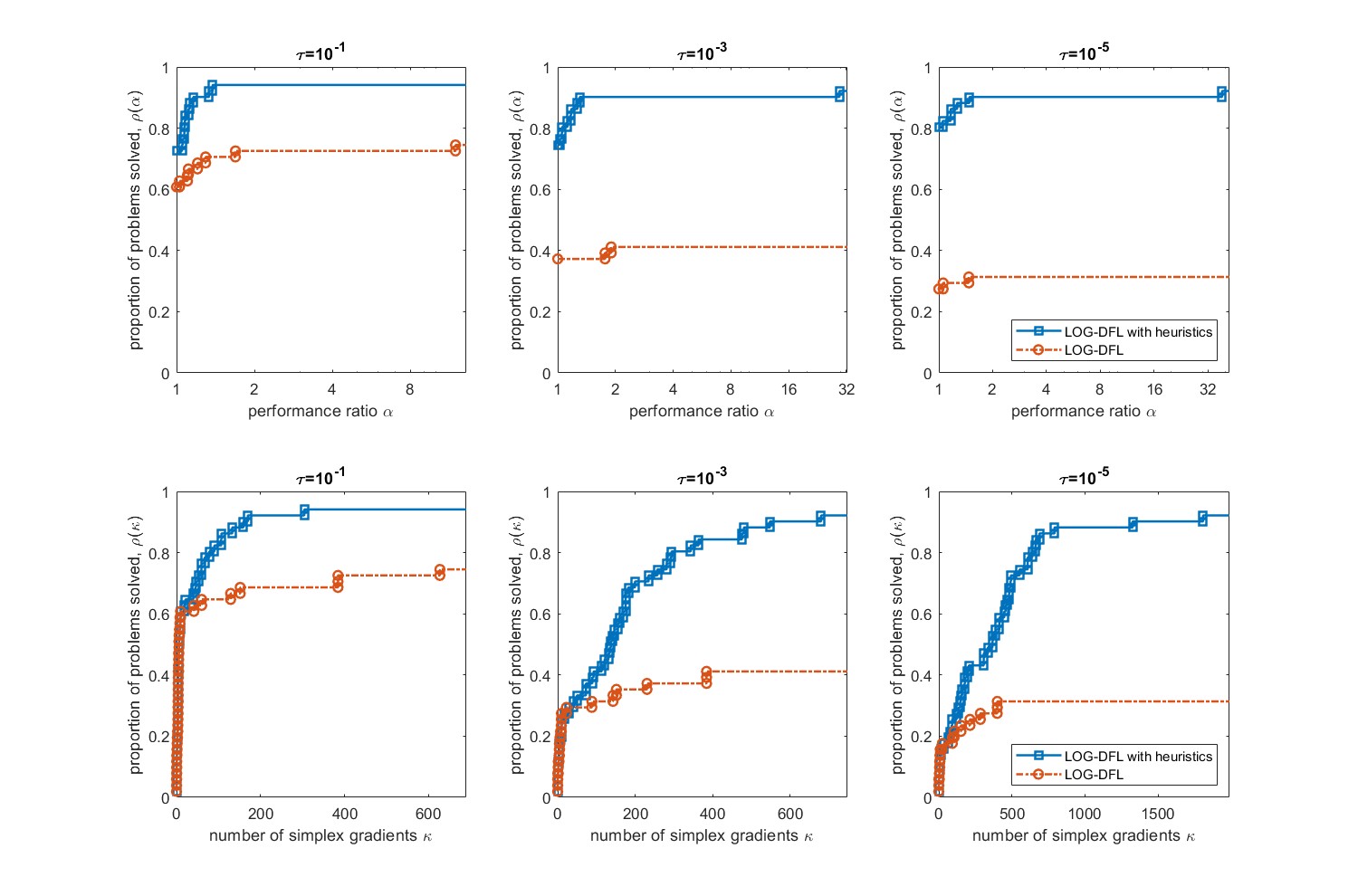}
        \caption{Problems where both find a feasible solution.}
        \label{fig:feastutti_innercompare}
    \end{subfigure}
    \caption{Performance and data profiles for the comparison between LOG-DFL with heuristic and LOG-DFL}
    \label{fig:innercompare}
\end{figure}

\noindent 
The comparison between our original algorithm LOG-DFL and the improved method (which we call LOG-DFL with heuristics) are reported in Figures \ref{fig:restutti_innercompare} and \ref{fig:feastutti_innercompare}. As we can see, the improved version of LOG-DFl is significantly better than the original version both in terms of efficiency and robustness.   

\begin{figure}[!ht]
    \centering
    \begin{subfigure}[b]{0.9\textwidth}
    \centering
        \includegraphics[width=0.8\textwidth]{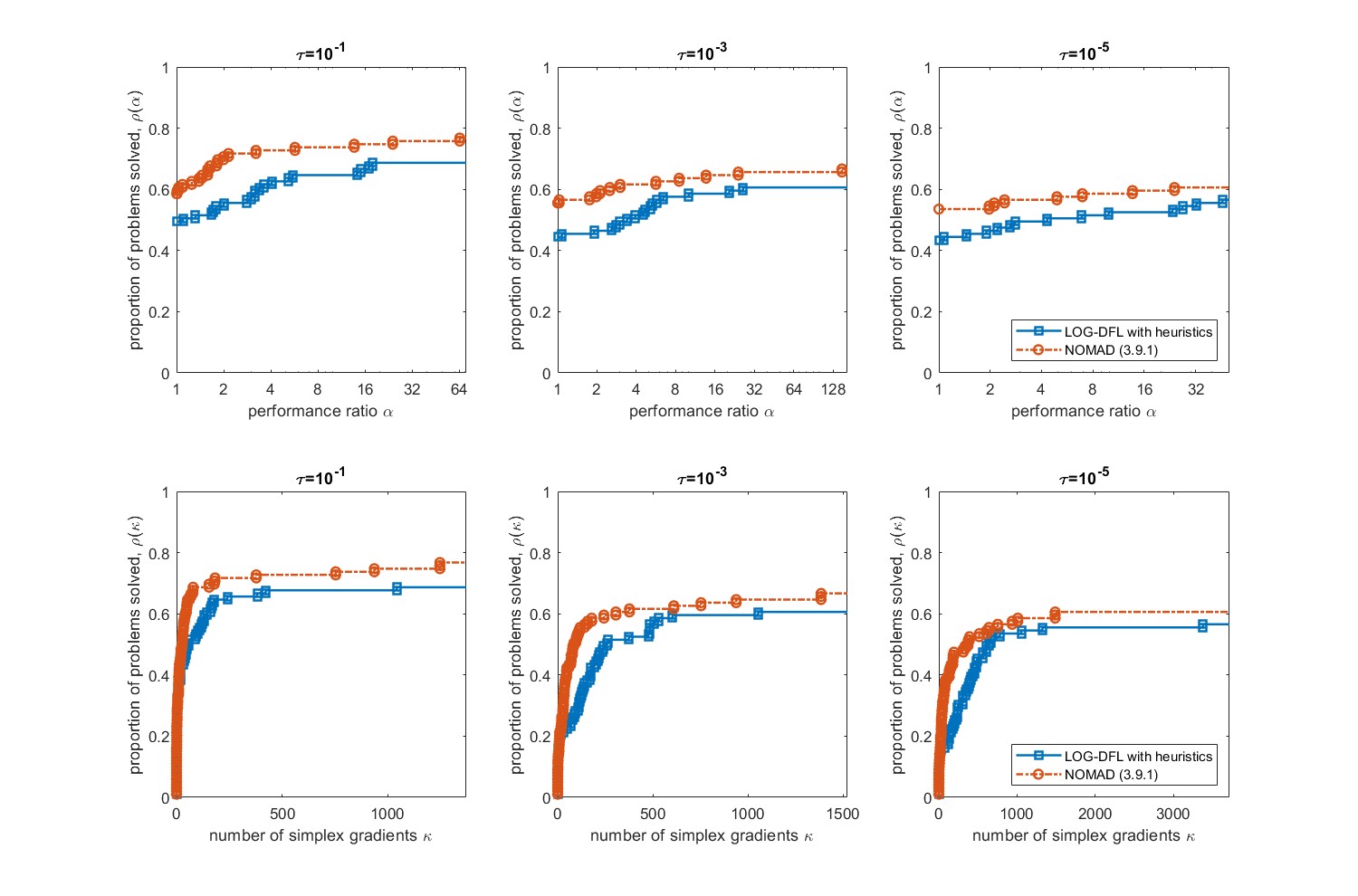}
        \caption{Results on the whole test problems collection.}
        \label{fig:restutti_bestNOMAD}
    \end{subfigure}
    \begin{subfigure}[b]{0.9\textwidth}
        \centering
        \includegraphics[width=0.8\textwidth]{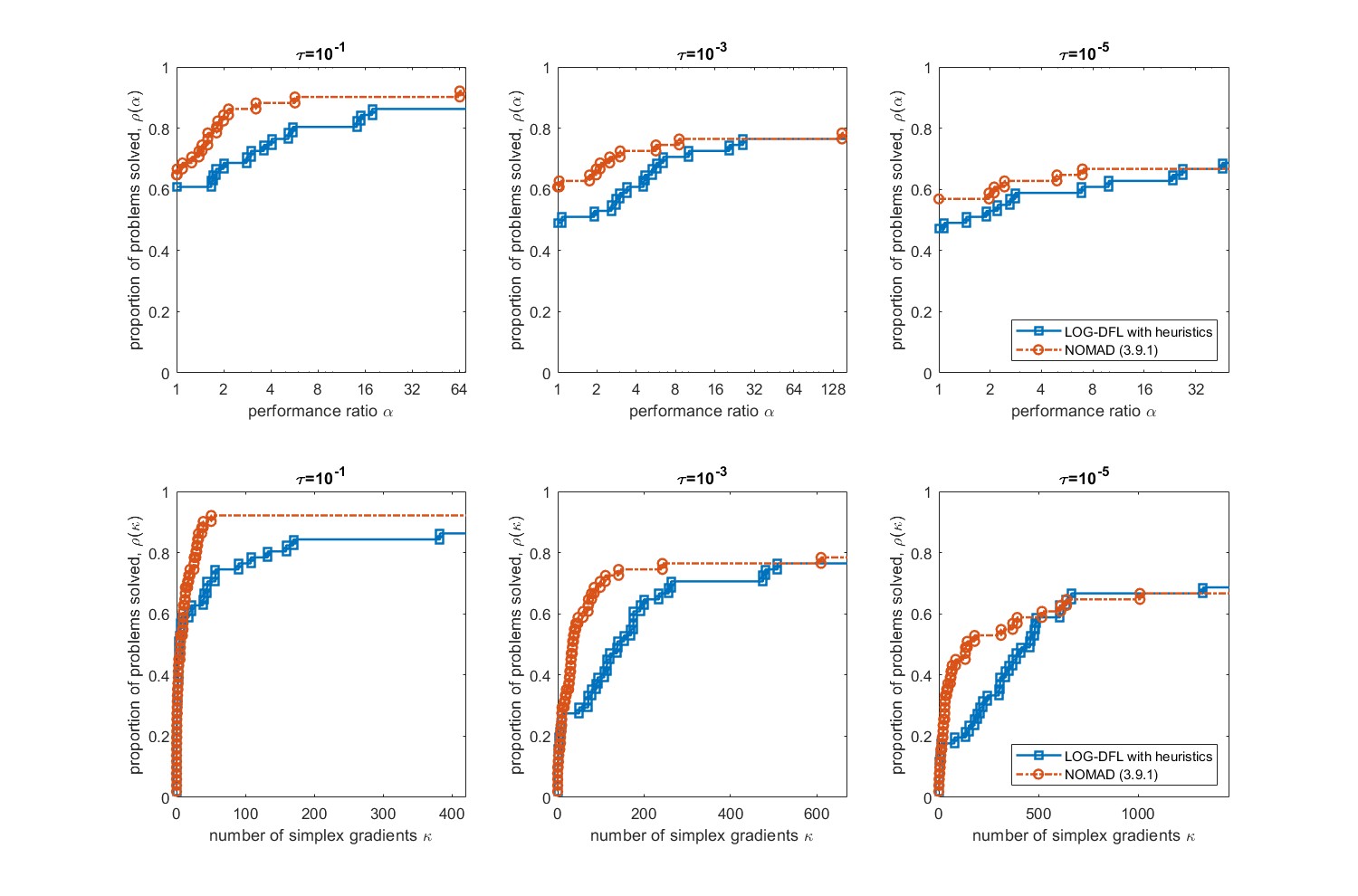}
        \caption{Results on problems where both methods find a feasible solution.}
        \label{fig:feastutti_bestNOMAD}
    \end{subfigure}
    \caption{Performance and data profiles for the comparison between LOG-DFL with heuristic and NOMAD}
    \label{fig:bestNOMAD}
\end{figure}

In Figures \ref{fig:restutti_bestNOMAD} and \ref{fig:feastutti_bestNOMAD}, we report the comparison between our best algorithm, i.e. LOG-DFL with heuristics, and NOMAD (3.9.1 with its default settings). As we can see, NOMAD outperforms LOG-DFL even though we are still comparable in terms of robustness. However, it must be noted that the default version of NOMAD has quadratic models enabled. Hence, the above comparison is not that fair after all. In fact, we are comparing our globalization strategy (based on the coordinate directions exploration) with the globalization strategy of NOMAD (i.e. the poll step) plus the search performed in the search step. This is something that could also be incorporated in our algorithm by suitably modifying and integrating the last step of LOG-DFL where the new point $x_{k+1}$ is defined. That said, it would be interesting to compare our method with NOMAD when models are disabled.   

\begin{figure}[!ht]
    \centering
     \begin{subfigure}[b]{0.9\textwidth}
         \centering
         \includegraphics[width=0.95\textwidth]{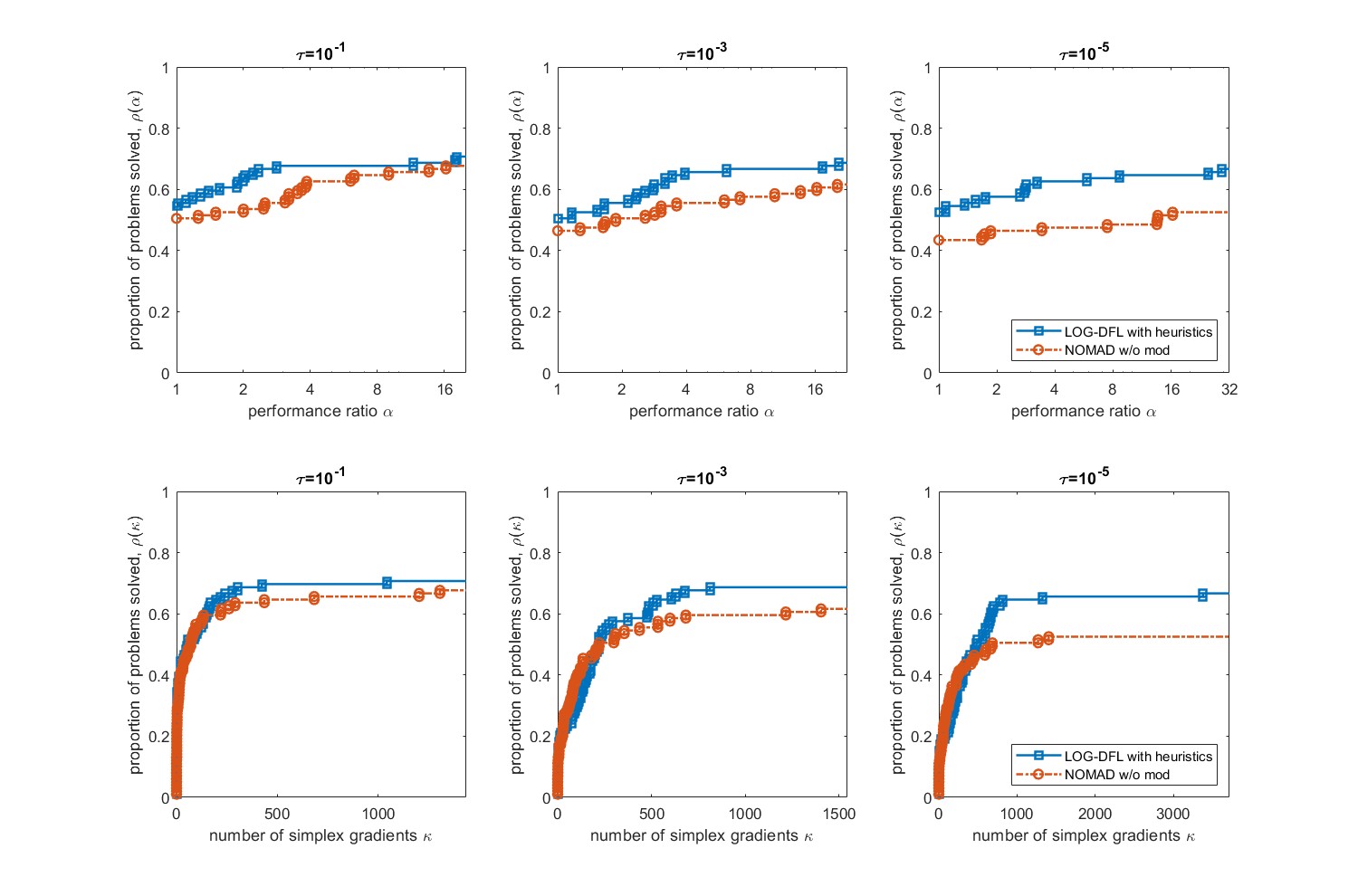}
         \caption{Whole test problems collection.}
         \label{fig:restutti_womod}
     \end{subfigure}
     \begin{subfigure}[b]{0.9\textwidth}
         \centering
         \includegraphics[width=0.95\textwidth]{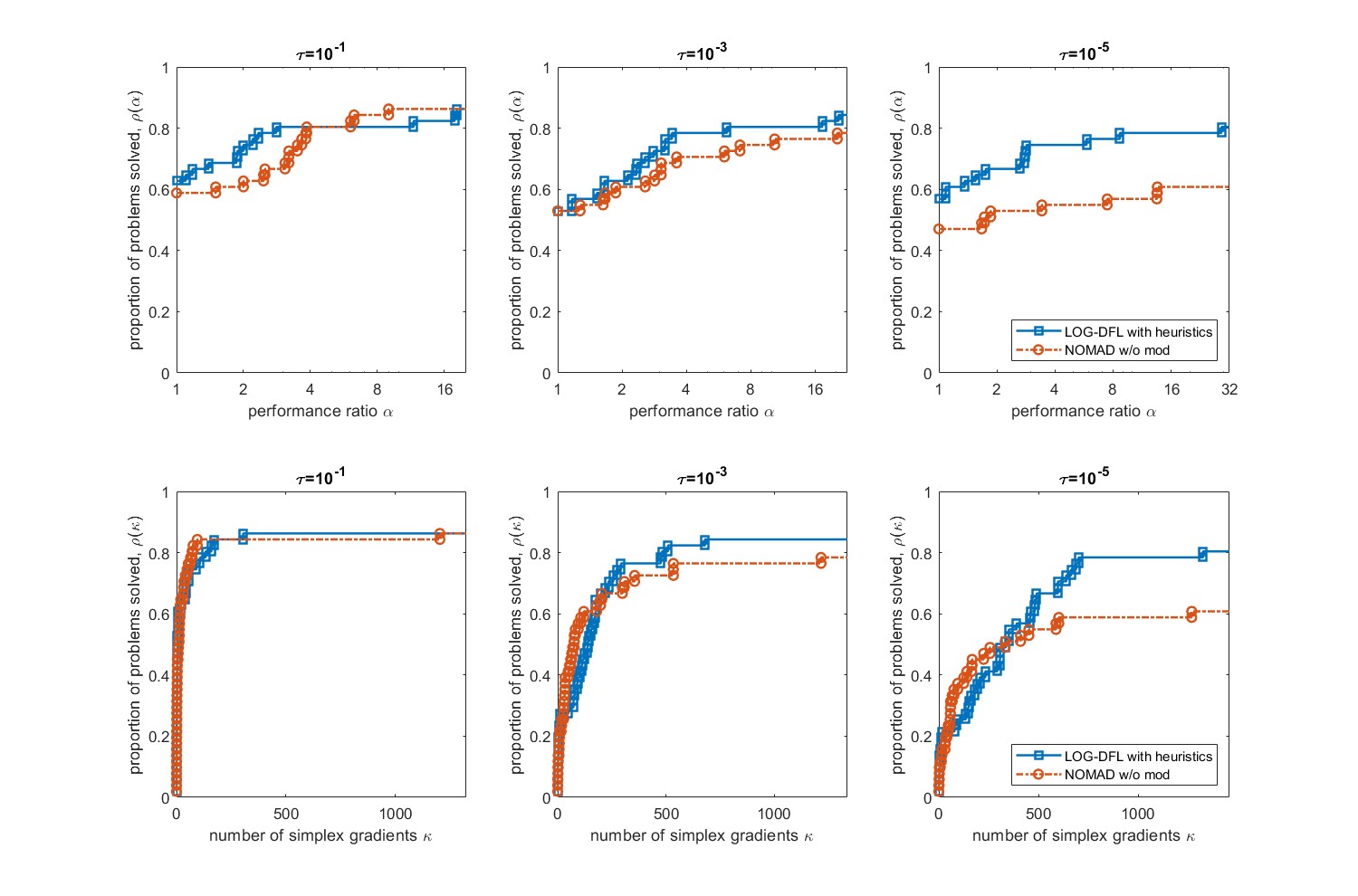}
         \caption{Problems where both find a feasible solution.}
         \label{fig:feastutti_womod}
    \end{subfigure}
    \caption{Performance and data profiles for the comparison between LOG-DFL with heuristic and NOMAD without models.}
    \label{fig:bestNOMAD1}
\end{figure}

Such a comparison is reported in Figures \ref{fig:restutti_womod} and \ref{fig:feastutti_womod}. It can be seen that when models are disabled NOMAD is outperformed by LOG-DFL (with heuristics). The superiority of LOG-DFL is more evident for high precision levels. 

%
%
%

\clearpage

\clearpage

\section{Concluding remarks}
\label{sec:conclusion}

In this paper we proposed a new algorithm based on the use of a mixed penalty-barrier merit function for the solution of constrained black-box problems. In particular, non-relaxable inequality constraints are handled by means of a log-barrier penalization. The main algorithm LOG-DFL is basically composed of three main steps.
\begin{enumerate}
    \item {\em Step 1} which is a modification of a quite standard search step for the minimization of the merit function. Note that the modification is indeed due to take into account the fact that (a subset of the) inequality constraints must be strictly satisfied.
    \item {\em Step 2} which is devoted to the barrier parameter updating criterion. The importance of this step emerges in the theoretical analysis and it was quite extensively described in Section \ref{sec:4}, just after the algorithm.
    \item {\em Step 3} which is the final step of the Algorithm. This is where the new iterate is chosen. It is worth noting that the new iterate can virtually be any point within $\stackrel{\circ}{\cal S}\cap X$ improving the merit function with respect to the point produced by Step 1. This particular step gives great freedom to our algorithm allowing for the use of any heuristic strategy that could take advantage of the points thus far produced (for instance, advanced modelling techniques could be used to approximate the merit function around the current point).
\end{enumerate}

\noindent For the proposed LOG-DFL Algorithm, we managed to prove convergence toward stationary points of the problem under quite mild assumptions. The convergence proof hinges on the crucial barrier parameter updating criterion (performed in Step 2 of the algorithm).

\par\noindent Furthermore, we also report a numerical experience and comparison with state-of-the-art solver on a large set of test problems from the CUTEst test set. The numerical results and comparison show that the proposed algorithm is both efficient and robust.
\par\medskip\noindent
{We note that the proposed algorithm and its theoretical properties can be easily adapted to optimization problems with more complex structures than \eqref{prob1}. In particular, inequality constraints violated at the starting point could be present and treated with the external penalization approach. Finally, the LOG-DFL algorithm is freely available for download through the DFL library as package {\tt LOGDFL} at the URL \href{http://www.iasi.cnr.it/~liuzzi/DFL/}{http://www.iasi.cnr.it/$\sim$liuzzi/DFL/} 
} 

\appendix

\section{Technical results}
\label{app:tech}
First we recall  a  result concerning a property of sequences of nonzero scalars which will be used in the proof of the
next proposition.

\begin{lemma}[see \cite{liuzzi2010sequential}]
\label{dgnew}
Let $\{a_k^{i}\}$, $i=1,\ldots ,p$, be sequences of nonzero scalars. There exist an index $i^\star\in \{1,\ldots ,p\}$ and an
infinite subset $K\subseteq \{0,1,\ldots \}$ such that
\begin{equation}
\label{confseq}
\lim_{k\to\infty,k\in K}\frac{a_k^{i}}{|a_k^{i^\star}|}=z_i,\quad\quad |z_i|<+\infty,\quad\quad i=1,\ldots ,p.
\end{equation}
\end{lemma}
\unskip





\par\medskip
Then, we report a technical result related to the behavior of Algorithm LOG-DFL which is necessary to prove boundedness of the {approximations of  multipliers introduced in Theorem~$\ref{mainalgdfl}$}.

\begin{proposition}
\label{theo3-lambda} Let the assumptions of
Theorem~$\ref{mainalgdfl}$ be satisfied {and let  $K$ be the set of indices defined in (\ref{definition_of_K}). If}
\[
\la_l(x;\eps) = -\frac{\eps}{g_l(x)},\qquad l=1,\dots,m.
\]
\[
\mu_j(x;\eps) = \frac{\nu}{\eps} \ \left|h_j(x)\right|^{\nu-1}, \qquad j=1,\dots,q
\]
then the subsequences $\{\la_l(x_k;\eps_k)\}_{K}$, $l=1,\dots,m$, and $\{\mu_j(x_k;\eps_k)\}_{K}$, $j=1,\dots,q$ are bounded.
\end{proposition}

\begin{proof}
By Propositions \ref{part_a_theorem}, we have that \eqref{d0bis_old}, \eqref{d1bis1_old}, \eqref{xi_ki_to_zero} and \eqref{y_ki_to_xbar} hold.
\par\medskip\noindent
{Let $\bar x$ be a limit point of the sequence $\{x_k\}_K$,
then there exists a subset of $K$,
which we relabel again $K$, such that
\[
\lim_{k\to\infty, k\in \bar K} x_k = \bar x
\]}
Let us denote $\bar D=D\cap D(\bar x)$. By applying the mean-value theorem to (\ref{new_fallimento}), we can write
$$
-o(\bar\xi_k^{i})\le P(y_k^i+\bar\xi_k^{i} d^{i};\eps_k) -
P(y_k^i;\eps_k)= \bar\xi_k^{i}\nabla
P(u_k^i;\eps_k)^{T}d^{i}\quad\quad \forall d^i\in \bar D,
$$
where $u_k^i=y_k^i+t_k^{i}\bar\xi_k^{i} d^{i}$, with $t_k^{i}\in (0,1)$. Thus, we have
\[
-\frac{o(\bar\xi_k^{i})}{\bar\xi_k^{i}} \leq \nabla P(u_k^i;\eps_k)^{T}
d^i\quad\quad \forall d^i\in \bar D.
\]
By considering the expression of $P(x;\eps)$, we can write
\begin{eqnarray}
\label{prop2.eq1} \lefteqn{\quad \nabla P(u_k^i;\eps_k)^T d^i
}\\\nonumber &&\quad = \Bigg(\nabla
f(u_k^i)+\sum_{l=1}^m\frac{\eps_k}{-g_l(u_k^i)}\nabla
g_l(u_k^i) \\\nonumber 
&& \qquad+ \sum_{j=1}^q \frac{\nu}{\eps_k} \left|h_j(u_k^i)\right|^{\nu-1}\nabla h_j(u_k^i) \Bigg)^{T} d^{i} \geq
-\frac{o(\bar\xi_k^{i})}{\bar\xi_k^{i}}\quad\quad \forall d^i\in \bar D.
\end{eqnarray}
Recalling that $u_k^i=y_k^i+t_k^{i}\xi_k^{i} d^{i}$, with $t_k^{i}\in (0,1)$, we have that
\begin{equation}
\label{prop2.eq5} \lim_{k\to\infty,k\in K} u_k^i = \bar x.
\end{equation}

By recalling the expression of $\la_l(x;\eps)$, $l=1,\dots,m$, and the expression of $\mu_j(x;\eps)$, $j=1,\dots,q$, we
can rewrite relation (\ref{prop2.eq1})~as
\begin{eqnarray}
\label{prop2.eq3} && \Bigg(\nabla f(u_k^i) + \sum_{l=1}^m\la_l(u_k^i;\eps_k)\nabla g_l(u_k^i) \\\nonumber
&& \qquad +\sum_{j=1}^q\mu_j(u_k^i;\eps_k)\nabla h_j(u_k^i) \Bigg)^T d^{i}
\geq -\frac{o(\bar\xi_k^{i})}{\bar\xi_k^{i}}\quad\quad \forall i: d^i\in
\bar D.
\end{eqnarray}
\par\medskip\noindent
First we prove that
\begin{equation}
\label{limite_lambda} \lim_{k\to\infty,k\in K}\left|\la_l(u_k^i;\eps_k)-
\la_l(x_k;\eps_k)\right|=0,\qquad l=1,\dots,m\quad \forall\
i:d^i\in \bar D.
\end{equation}
In fact,
\begin{eqnarray}
\label{mind4}&& \qquad
\left|\frac{\eps_k}{-g_l(u_k^i)}-
\frac{\eps_k}{-g_l(x_k)}\right|=
\eps_k\left|\frac{g_l(x_k)-g_l(u_k^i)}{(-g_l(u_k^i))(-g_l(x_k))}\right|=
\\\nonumber
 &&\qquad  =\eps_k\frac{\left|\nabla g_l(\tilde u_k^{i,l})^T(x_k - u_k^i) \right|}{g_l(u_k^i)g_l(x_k)}\leq \eps_k\frac{\|\nabla g_l(\tilde u_k^{i,l})\|\|(u_k^i-x_k)\|}{|g_l(u_k^i)||g_l(x_k)|},
\end{eqnarray}
where $\tilde u_k^{i,l} = u_k^i + \tilde t_k^{i,l}x_k$ with $\tilde t_k^{i,l}\in (0,1)$. Then,
\begin{eqnarray}
\nonumber
&& \eps_k\frac{\|\nabla g_l(\tilde u_k^{i,l})\|\|(u_k^i-x_k)\|}{|g_l(u_k^i)||g_l(x_k)|} \leq c_1 \eps_k\frac{\max_{i:d^i\in\bar D}\{\bar\xi_k^i,\|y_k^i-x_k\|\}}{|g_l(u_k^i)||g_l(x_k)|}
\end{eqnarray}
Now, we show that that 
\begin{equation}
\label{appendice_utile1}\frac{1}{|g_l(u_k^i)|} \le c_2\frac{1}{|g(y_k^i)|}, \quad \forall k, k\in K.
\end{equation}
To this end, we assume by contradiction that there does not exist such a constant $c_2$.
This would imply:
\begin{eqnarray}
\label{absurd_u}
\lim_{k\to\infty,k\in K}\frac{\frac{1}{|g_l(u_k^i)|}}{\frac{1}{|g_l(y_k^i)|}} = \lim_{k\to\infty, k \in K}\frac{|g_l(y_k^i)|}{|g_l(u_k^i)|} = +\infty,
\end{eqnarray}
and let us consider the case:
\begin{eqnarray}
\label{g_to_0}
\lim_{k\to\infty,k\in K} |g_l(y_k^i)|= 0.
\end{eqnarray}
Since $g_l(y_k^i)<0$ and $g_l(u_k^i)<0$ $\forall k$, by (\ref{absurd_u}) $\bar{k}$ exists such that, for all $k\ge\bar{k}$, $k\in K$, we have:
\[
-g_l(y_k^i) > -g_l(u_k^i).
\]
Recalling the definition of $u_k^i$:
\begin{eqnarray}
\label{absurd_y}
-g_l(y_k^i) > -g_l(y_k^i+t_k^i \bar\xi_k^i d^i).
\end{eqnarray}
Using the Lipschitz continuity assumption on $g$, $||d^i||=1$, (\ref{absurd_y}) and recalling the possible choices of $\bar\xi_k^i$ described in Proposition \ref{part_a_theorem}:
\begin{eqnarray}
\label{absurd_y_replace}
&&-g(y_k^i + t_k^i \bar\xi_k^i d^i) > -g(y_k^i) - Lt_k\bar\xi_k^i \ge\\ \nonumber
&& \ge  -g(y_k^i) - Lt_k\max_{i=1,2,\dots,n}\{\tilde\alpha_k^i,\alpha_k^i\} , \quad \forall k\ge\bar{k}, k\in K.
\end{eqnarray}
The instructions of Step 2 imply that, for all $k\in K$:
\[
\max_{i=1,2,\dots,n}\{\tilde\alpha_k^i,\alpha_k^i\}\le\min\{
 \epsilon_k^p,  (g_{\min})_k^2\}.
\] 
Hence, $\bar{k}_2$ exists such that:
\begin{eqnarray}
\label{positivity_g_alfa}
&&-g(y_k^i) - Lt_k\max_{i=1,2,\dots,n}\{\tilde\alpha_k^i,\alpha_k^i\}\ge\\ \nonumber
&&\ge -g(y_k^i) - Lt_k(g(y_k^i)^2) \ge 0 , \quad \forall k\ge\max\{\bar{k},\bar{k}_2\}, k\in K.
\end{eqnarray}
That allows us to say:
\begin{eqnarray}\nonumber
&&\lim_{k\to\infty, k \in K}\frac{-g_l(y_k^i)}{-g_l(u_k^i)} = \lim_{k\to\infty}\frac{-g_l(y_k^i)}{-g_l(y_k^i+t_k^i\xi_k^id^i)}\\ \nonumber
&& \le \lim_{k\to\infty, k \in K}\frac{-g(y_k^i)}{-g(y_k^i) - Lt_k(g(y_k^i)^2)} = 1,
\end{eqnarray}
which leads to a contradiction thus proving (\ref{appendice_utile1}).\\
Let us now consider the case:
\[
\lim_{k\to\infty,k\in K} |g_l(y_k^i)| = c < +\infty,
\]
which implies:
\begin{eqnarray}\nonumber
\lim_{k\to\infty, k \in K}\frac{-g_l(y_k^i)}{-g_l(u_k^i)} = \lim_{k\to\infty}\frac{-g_l(y_k^i)}{-g_l(y_k^i+t_k^i\xi_k^id^i)} < \infty,
\end{eqnarray}
which again leads to a contradiction proving (\ref{appendice_utile1}).\\
\par\medskip\noindent
Hence, the existence of the constant $c_2$, (\ref{appendice_utile1}), and recalling that $\bar\xi_k^i\leq\xi_k^i$, allow us to write
\begin{eqnarray}
\nonumber
&& c_1 \eps_k\frac{\max_{i:d^i\in\bar D}\{\bar\xi_k^i,\|y_k^i-x_k\|\}}{|g_l(u_k^i)||g_l(x_k)|} \le c_1 \eps_k\frac{\max_{i:d^i\in\bar D}\{\xi_k^i,\|y_k^i-x_k\|\}}{c_2|g_l(y_k^i)||g_l(x_k)|}\\ \nonumber
&&  \le c_1 \eps_k\frac{\max_{i:d^i\in\bar D}\{\xi_k^i,\|y_k^i-x_k\|\}}{c_2\min_{i:d^i\in\bar D, l}\{|g_l(y_k^i)|,|g_l(x_k)|\}^2}.
\end{eqnarray}
\par\medskip\noindent
Now, recalling 
that $y_k^i = x_k + \sum_{j=1}^{i-1}\alpha_k^j d_k^j$ {and the possible choices for $\xi_k^i$ described in  Proposition \ref{part_a_theorem}} and the definition of $\bar\xi_k^i$,
we can write
\begin{eqnarray}
\label{appoggio2}
&&{\eps_k\frac{\max_{i:d^i\in\bar D}\{\xi_k^i,\|y_k^i-x_k\|\}}{\min_{i:d^i\in\bar D, j}\{|g_j(y_k^i)|,|g_j(x_k)|\}^2}} \\\nonumber
&& \qquad \le n \eps_k \frac{\max_{i=1,2,\dots,n}\{\tilde\alpha_k^i,\alpha_k^i\}}{(g_{min})_k^2}.
\end{eqnarray}
\par\medskip\noindent
The instructions of Step 2 imply
that, for all $k\in K$,
\[
\max_{i=1,2,\dots,n}\{\tilde\alpha_k^i,\alpha_k^i\}\le\min\{
 \epsilon_k^p,  (g_{\min})_k^2\}
\] 
so that
\[
 n \eps_k \frac{\max_{i=1,2,\dots,n}\{\tilde\alpha_k^i,\alpha_k^i\}}{(g_{min})_k^2} \leq n\eps_k
\]
Then, (\ref{limite_lambda}) is proved by (\ref{mind4}) and recalling Proposition \ref{eps_to_0}.
\par\medskip\noindent
Now, we prove that:
\begin{equation}
\label{limite_mu} \lim_{k\to\infty,k\in K}\left|\mu_j(u_k^i;\eps_k)-
\mu_j(x_k;\eps_k)\right|=0,\qquad j=1,\dots,q\quad \forall\
i:d^i\in \bar D.
\end{equation}
In fact,
\begin{eqnarray}
\label{mind5}&& \qquad
\left|\; \Bigg|\frac{h_j(u_k^i)}{\eps_k}\Bigg|^{\nu-1} -\; \Bigg|\frac{h_j(x_k)}{\eps_k}\Bigg|^{\nu-1}\right|
\\[1em]\nonumber
 &&\qquad  = \left| \; \Bigg|\frac{h_j(x_k)}{\eps_k}+\frac{1}{\eps_k}\nabla h_j(\tilde u_k^{i,l})^T(u_k^i - x_k)\Bigg|^{\nu-1} -\; \Bigg|\frac{h_j(x_k)}{\eps_k}\Bigg|^{\nu-1}\right| 
\\[1em]\nonumber
&&\qquad  \le \left|  \Bigg|\frac{h_j(x_k)}{\eps_k}\ \Bigg|^{\nu-1}+ \Bigg|\frac{1}{\eps_k}\nabla h_j(\tilde u_k^{i,l})^T(u_k^i - x_k) \Bigg|^{\nu-1} - \;  \Bigg|\frac{h_j(x_k)}{\eps_k}\Bigg|^{\nu-1}\right|
\\[1em]\nonumber
&&\qquad  = \Bigg|\frac{1}{\eps_k}\nabla h_j(\tilde u_k^{i,l})^T(u_k^i - x_k)\Bigg|^{\nu-1} \le \frac{\|\nabla h_j(\tilde u_k^{i,l})\|^{\nu-1}\|(u_k^i-x_k)\|^{\nu-1}}{e_k^{\nu-1}},
\end{eqnarray}
where again $\tilde u_k^{i,l} = u_k^i + \tilde t_k^{i,l}x_k$ with $\tilde t_k^{i,l}\in (0,1)$. Now, recalling that $h_i$, $i=1,\dots,q$, are continuously differentiable functions and that $u_k^i=y_k^i+t_k^{i}\bar\xi_k^{i} d^{i}$, with $t_k^{i}\in (0,1)$, from \eqref{mind5} and $\bar\xi_k^i\leq\xi_k^i$ we can write
\begin{equation}
\label{mind5_1}
\left|\; \Bigg|\frac{h_j(u_k^i)}{\eps_k}\Bigg|^{\nu-1} -\; \Bigg|\frac{h_j(x_k)}{\eps_k}\Bigg|^{\nu-1}\right| \le c_2 \left(\frac{\max_{i: d^i \in \bar D}\{\xi_k^i,\|y_k^i - x_k\|\}}{\eps_k}\right)^{\nu-1}
\end{equation}

\par\noindent
Note that
\[
\frac{\max_{i: d^i \in \bar D}\{\xi_k^i,\|y_k^i - x_k\|\}}{\eps_k} \leq \frac{\max_{i=1,2,\dots,n}\{\tilde\alpha_k^i,\alpha_k^i\}}{\eps_k}.
\]
By the instructions of Algorithm LOG-DFL, for $k\in K$ we can write
\[
\max_{i=1,2,\dots,n}\{\tilde\alpha_k^i,\alpha_k^i\} \leq \min\{\eps_k^p,(g_{\min})_k^2\}
\]
that is
\[
\frac{\max_{i=1,2,\dots,n}\{\tilde\alpha_k^i,\alpha_k^i\}}{\eps_k} \leq \eps_k^{p-1}.
\]
Then, (\ref{limite_mu}) is proved by (\ref{mind5_1}) and recalling Proposition \ref{eps_to_0}.

\par\medskip\noindent
Now, we are ready to show that the sequences $\{\la_l(x_k;\eps_k)\}_{K}$, $l=1,\dots,m$, and $\{\mu_j(x_k;\eps_k\}_{K}$, $j=1,\dots,p$ are bounded.

\noindent In fact, by simple manipulations (\ref{prop2.eq3}) can be rewritten as
\begin{eqnarray}
\label{prop2.eq3-1}
\lefteqn{\Bigg( \nabla f(u_k^i) + \sum_{l=1}^m\nabla g_l(u_k^i)\la_l(x_k;\eps_k) } \\\nonumber
 & + & \sum_{l=1}^m\nabla
g_l(u_k^i)\left(\la_l(u_k^i;\eps_k)-\la_l(x_k;\eps_k)\right)+\sum_{j=1}^q \nabla h_j(u_k^i)\mu_j(x_k;\eps_k) \\\nonumber
& + & \sum_{j=1}^1\nabla
h_j(u_k^i)\left(\mu_j(u_k^i;\eps_k)-\mu_j(x_k;\eps_k)\right) \Bigg)^T d^{i} \geq -\frac{o(\bar\xi_k^{i})}{\bar\xi_k^{i}}\quad\quad
\forall i: d^i\in \bar D.
\end{eqnarray}
Let
$$
\{a_k^1,\ldots ,a_k^{m}\}=\{\lambda_1(x_k;\epsilon_k),\ldots ,\lambda_m(x_k;\epsilon_k)\}.
$$
$$
\{a_k^{m+1},\dots,a_k^{m+q}\}=\{\mu_1(x_k;\eps_k),\dots,\mu_q(x_k;\eps_k)\}
$$
By contradiction let us assume that there exists at least an index $h\in\{1,\ldots ,{m+}q\}$ such that
$$
 \lim_{k\to\infty,k\in K}|a_k^h|=+\infty.
$$
From Lemma~\ref{dgnew}, we get that there exist an infinite subset
(which we again relabel  $K$) and an index $s\in\{1,\ldots ,{m+}q\}$ such
that,
\begin{equation}
\label{l1}
\lim_{k\to\infty,k\in K}{\frac{ a_k^i}{|a_k^s|}}=z_i,\quad |z_i| < +\infty,\quad i=1,\dots,q.
\end{equation}
Note that
\begin{equation}
\label{seqinf1}
{z_i \ge 0,\quad i\in\{1,\ldots ,m\},} \quad\quad z_s =1\quad\quad {\rm and}\quad\quad |a_k^s|\to +\infty.
\end{equation}
Dividing relation (\ref{prop2.eq3-1}) by $|a_k^s|$, we have
\begin{eqnarray}
\label{prop2.eq3-2}
\lefteqn{\Bigg(\frac{\nabla f(u_k^i)}{|a_k^s|} + \sum_{l=1}^m\frac{\nabla g_l(u_k^i)a_k^l}{|a_k^s|} } \\\nonumber
 &\quad + & \sum_{l=1}^m\nabla
g_l(u_k^i)\frac{\la_l(u_k^i;\eps_k)-\la_l(x_k;\eps_k)}{|a_k^s|}
+\sum_{j=1}^q \frac{\nabla h_j(u_k^i)a_k^{m+j}}{|a_k^s|} \\\nonumber
& + & \sum_{j=1}^1\nabla
h_j(u_k^i)\frac{\mu_j(u_k^i;\eps_k)-\mu_j(x_k;\eps_k)}{|a_k^s|} \Bigg)^T d^{i} \geq -\frac{o(\bar\xi_k^{i})}{\bar\xi_k^{i}}\quad\quad
\forall i: d^i\in \bar D.
\end{eqnarray}
Taking the limits for $k\to\infty$ and $k\in K$, recalling that
$|a_k^s|\to\infty$, and using (\ref{limite_lambda}), (\ref{limite_mu}), (\ref{l1}), and
(\ref{prop2.eq5}), we obtain
\begin{equation}
\label{p1} \left(\sum_{l=1}^mz_l\nabla g_l(\bar x) + \sum_{j=1}^q z_{m+j}\nabla h_j(\bar x)\right)^Td^i \ge 0\quad \quad
\forall i:d^i\in \bar D.
\end{equation}
We recall that, $\bar x$ satisfies the MFCQ by assumption. Now, let
$\hat d\in D(\bar x)$ be the direction considered in
Definition~\ref{ass}, which, 
from Proposition~\ref{cn1}, can be written as
\begin{equation}
\label{combcone1}
\hat d = \sum_{i:d^i\in \bar{D}}\hat \beta_i d^i.
\end{equation}
Thus, from (\ref{combcone1}) and (\ref{p1}), we obtain
{
\begin{eqnarray}
\label{dirassnew1}
\lefteqn{\left(\sum_{l=1}^mz_l\nabla
g_l(\bar x)+ \sum_{j=1}^q z_{m+j}\nabla h_j(\bar x) \right)^T\hat d \ =} \\\nonumber
&&  = \sum_{l=1}^mz_l\nabla g_l(\bar x)^T\hat d + \sum_{j=1}^q z_{m+j}\nabla h_j(\bar x)^T\hat d \ge 0.
\end{eqnarray}
The above relation, considering definition~\ref{ass}, implies
\begin{eqnarray}
\label{dirassnew2}
 \sum_{l=1}^mz_l\nabla g_l(\bar x)^T\hat d \ge 0.
\end{eqnarray}
By definition, we note that 
\begin{eqnarray}
\label{dirassnew3}
&& z_i=0  \qquad \hbox{for all}\quad i \in \{1,\ldots,m\} \quad \hbox{and} \quad i\notin I^+(\bar x).
\end{eqnarray}
Furthermore, by Definition~\ref{ass}, (\ref{dirassnew2}) and (\ref{dirassnew3}) imply 
\begin{eqnarray}
\label{dirassnew4}
&& z_i=0  \qquad \hbox{for all}\quad i \in \{1,\ldots,m\} \quad \hbox{and} \quad i\in I^+(\bar x).
\end{eqnarray}
Hence, recalling \eqref{p1}, \eqref{dirassnew3} and \eqref{dirassnew4} we have that
\begin{equation}
\label{p2new} 
\left( \sum_{j=1}^q z_{m+j}\nabla h_j(\bar x)\right)^Td^i \ge 0\quad \quad
\forall i:d^i\in \bar D.
\end{equation}
By using again Definition~\ref{ass},  Proposition~\ref{cn1} and \eqref{p2new}, we obtain
\begin{eqnarray}
\label{dirassnew5}
&& z_{m+j}=0  \qquad \hbox{for all}\quad j \in \{1,\ldots,q\}.
\end{eqnarray}
In conclusion we get that \eqref{dirassnew2}, \eqref{dirassnew3} and
\eqref{dirassnew5} contradict \eqref{seqinf1}.
and this concludes the proof.}
\end{proof}


\end{document}